%% file: TCST-2017-0381_R1_MX_PKP_DC_20170908_MX_no_note.tex
\documentclass[journal]{IEEEtran}
\ifCLASSINFOpdf
\else
\fi
\input{sy-headinput}
\usepackage{amsmath}
\usepackage{graphicx}
\usepackage{amssymb}
\usepackage{cite}
\usepackage{cases}

\usepackage{multicol}
\usepackage[figuresright]{rotfloat}[2004/01/04]
\usepackage{amsfonts}
\usepackage{setspace}
\usepackage{cite}
\usepackage{mathrsfs}
\usepackage{graphics,dblfloatfix} 
\usepackage{epsfig} 
\usepackage{mathdots}
\usepackage[usenames, dvipsnames]{color}
\allowdisplaybreaks[1]
\usepackage{algorithm}
\usepackage{algpseudocode}
\usepackage{soul}
\usepackage{xargs}                      

\begin{document}

%
\title{Decentralized Charging Control of Electric Vehicles in Residential Distribution Networks}
%
%
%

\author{Mingxi~Liu\IEEEauthorrefmark{2},~\IEEEmembership{Member,~IEEE,}
        Phillippe~K.~Phanivong\IEEEauthorrefmark{2},~\IEEEmembership{Student~Member,~IEEE,}
        Yang~Shi\IEEEauthorrefmark{4},~\IEEEmembership{Fellow,~IEEE,}
        Duncan~S.~Callaway\IEEEauthorrefmark{2},~\IEEEmembership{Member,~IEEE}
\thanks{\IEEEauthorrefmark{2}M. Liu, P. Phanivong, and D. Callaway are with the Energy \& Resources Group at University of California, Berkeley, 310 Barrows Hall, Berkeley, CA, 94720, USA. {\tt\small \{mxliu, phillippe\_phanivong, dcal\}@berkeley.edu}. }
\thanks{\IEEEauthorrefmark{4}Y. Shi is with the Department of Mechanical Engineering and Institute for Integrated Energy Systems at the University of Victoria, Victoria, 3800 Finnerty Road, BC, Canada V8W 2Y2.  {\tt\small yshi@uvic.ca}.}
\thanks{This work was supported by the Natural Sciences and Engineering Research Council of Canada, the California Energy Commission Award EPIC 15-013, the National Science Foundation CyberSEES No. 1539585, and the State Key Laboratory of Alternate Electrical Power System with Renewable Energy Sources (Grant No. LAPS17006).}
}

%



\maketitle

\begin{abstract}
Electric vehicle (EV) charging can negatively impact electric distribution networks by exceeding equipment thermal ratings and causing voltages to drop below standard ranges. In this paper, we develop a decentralized EV charging control scheme to achieve ``valley-filling'' (i.e., flattening demand profile during overnight charging), meanwhile meeting heterogeneous individual charging requirements and satisfying distribution network constraints. The formulated problem is an optimization problem with a non-separable objective function and strongly coupled inequality constraints. We propose a novel shrunken-primal-dual subgradient (SPDS) algorithm to support the decentralized control scheme, derive conditions guaranteeing its convergence, and verify its efficacy and convergence with a representative distribution network model.

\end{abstract}


%
\IEEEpeerreviewmaketitle

\section{Introduction} \label{Introduction}

\IEEEPARstart{C}{ombustion} of fossil-based fuel is one of the main contributors to the greenhouse gas (GHG) emissions in most countries. For example, in the US, the most recent data (2014) shows that 26.5\% of total GHG emissions comes from the transportation sector \cite{EPA_GHG_2016}. Electric vehicles (EVs) have been proved to be effective in increasing energy conversion efficiency and reducing GHG emissions by relieving the usage of fossil-based fuel \cite{Liu_PIEEE_2013, Chan_PIEEE_2007, Chau_PIEEE_2007}. From 2013 to 2015, the market share of EVs presented a significant growth in most countries \cite{Cobb_url_2014, Shahan_url_2013}. Provided they are charged with low-emissions power sources, the increasing number of EVs will no-doubt benefit the global environment.

Beyond GHG reduction, there are extensive studies on actively utilizing EVs for grid-level benefits. In \cite{Callaway_PIEEE_2011, Hiskens_VPPC_2009}, the authors envisioned the potential of EV charging control for the provision of grid services. Some studies have focused on modeling and control problems of EVs for valley-filling \cite{Gan_TPS_2013, Liu_PIEEE_2016,Ma_Automatica_2016, Zhan_EPSR_2015,Karfopoulos_TPS_2016}, load balancing \cite{Mercurio_MCCA_2013}, and frequency regulation \cite{Lam_TSG_2016, Karfopoulos_TPS_2016_2,Yang_TPS_2013, Liu_TPS_2015, Lin_IOT_2014}, to name a few. In the presence of renewable energy, authors in \cite{Richardson_RSER_2013, Kou_TSG_2016, Bashash_TSG_2012} discussed the potential of EV charging control for facilitating the integration of wind power and solar power.  Jointly controlling EVs, solar energy, and energy storage systems to achieve valley-filling, peaking shaving, and energy neutral design was studied in \cite{Denholm_JPS_2013,Munkhammar_APEN_2015,Wi_TCE_2013}. 

Distribution network considerations can be extremely important in this context, and the aforementioned papers neglect network physics and constraints. For example, if the charging process of a large population of EVs is not properly controlled, it would negatively affect the power system, such as the elevated existing demand peak and newly created demand peak \cite{Gan_TPS_2013}. Lopes \emph{et al.} \cite{Lopes_PIEEE_2011} elaborated on the challenges of EV integration into the traditional mid- or low-voltage distribution network, including severe nodal voltage drop, transformer overloading, network congestion, and increased power loss.  In \cite{Fernandez_TPS_2011}, through analyzing simulation results in different application scenarios, the authors verified that uncontrolled charging or controlled charging without considering network constraints would negatively affect the distribution network. Other literature reviewing the impacts includes \cite{Habib_JPS_2015, Clement_Nyns_TPS_2010}. 

Researchers have begun to develop methods to reduce the impact of EV charging on distribution networks.  In \cite{Geth_ISGT_2012}, the authors developed a voltage droop charging control to maintain the nodal voltage level. Quir$\acute{\text{o}}$s-Tort$\acute{\text{o}}$s \emph{et al.} \cite{Quiros-Tortos_TPS_2016} studied a centralized algorithm, currently being trialed in the UK, to mitigate voltage drop and transformer overloading in the low-voltage distribution network. In \cite{Bansal_CDC_2014}, the authors developed a centralized model predictive control (MPC) scheme to maintain the voltage profile while satisfying charging requirements. Although the above-mentioned results can help alleviate the impacts, they were designed for meeting the network constraints only. We are aware of only a limited number of papers that address the problem of distribution network-aware provision of grid services.   Richardson \emph{et al.} \cite{Richardson_TPS_2012} optimized the EV charging profiles to minimize the total power consumption. Luo and Chan \cite{Luo_IET_2013} studied a real-time control design based on the voltage profile leveling to minimize power losses. 

Among all types of grid services, a number of papers focus on filling the overnight load ``valley'', where the non-EV electricity is at its lowest \cite{Ma_TCST_2013, Kundu_CDC_2012}. By filling the valley, on the one hand, EVs could be charged during the night without causing any inconvenience; on the other hand, the daily operations of power plants and associated cost could be reduced \cite{Denholm_report_2006}. In our paper, we aim to develop a framework for the provision of valley-filling while satisfying both local charging needs and distribution network constraints. This framework enjoys the security and privacy features, the guaranteed service performance, low communication requirements, and can be readily extended to facilitate other grid services, such as power trajectory tracking. 

The core of establishing such a framework lies in the control scheme design. Most literature addressing charging control problems under network constraints utilized the centralized control scheme \cite{Clement_Nyns_TPS_2010, Richardson_TPS_2012, Sharma_TSG_2014, Hu_TSG_2014, Luo_IET_2013, Quiros-Tortos_TPS_2016, Bansal_CDC_2014}.  Though centralized schemes are easy to realize in algorithm design, they do not scale well in a computational sense.   In contrast, decentralized control distributes the heavy computing load to individual agents. Each agent only needs to solve its own problem of small size without communicating with others. It's worth mentioning that, when implementing the charging control in the real world, two competing sets of standards are defining the ways EVs can communicate in the grid: ISO/IEC 15118 \cite{ISO_15118}, a protocol developed mainly in Europe, includes a full billing system that communicates through the electric vehicle supply equipment (EVSE) to the EV; while in the US, SAE 2847 \cite{SAE_2847} was developed to integrate EVs into the SEP 2.0 protocol \cite{SEP_20}. Although both of these standards have different approaches to communication and control, a decentralized control scheme will reduce communication requirements and allow for more efficient computation. Additionally, both of these standards require the EV state-of-charge (SOC) information retained locally, which can only be achieved by a decentralized scheme. Hence, our objective in this paper is to design a decentralized optimal controller that (1) can be embedded into charging points/EVs, (2) does not share local SOC information with the utility or other EVs, and (3) does not require a communication network between EVs.

General decentralized/distributed optimization algorithms have been considered widely in the literature \cite{Jia_ACC_2001, Liu_TIE_2016,Rivera_CDC_2013,Chen_arXiv_2016, Gao_arXiv_2016, Cui_arXiv_2015,Li_TAC_2014}, however none of them can solve the valley-filling problem considered here: On the one hand, the valley-filling objective function is a coupled and non-separable one; on the other hand, nodal voltages strongly couple the individual charging power via an inequality constraint. Only a small number of papers have addressed the distributed/decentralized realization of this problem. Chang \emph{et al.} \cite{Chang_TAC_2014} proposed a consensus-based primal-dual perturbation algorithm for a distributed consensus problem. Koshal \emph{et al.} \cite {Koshal_SIAM_2011} developed a regularized primal-dual subgradient (RPDS) algorithm via regularizing both primal and dual variables in the Lagrangian. This RPDS guarantees convergency, however introduces relative errors to the optimal solution. Authors in \cite{Zhang_TPS_2017} proposed an Alternating Direction Method of Multipliers (ADMM)-based decentralized control scheme to tackle the same problem. This design is fundamentally applicable as it enjoys privacy and security features, guarantees convergence, and poses the minimal computational requirements to EV controllers. However, the required two-layer communication network complicates the communication and poses computing burdens to all buses. In addition, ADMM-based decentralized control schemes always face the problem of a large number of iterations. Motivated by the above facts, we solve this decentralized charging control problem by proposing a novel shrunken-primal-dual subgradient (SPDS) algorithm, which eliminates unnecessary convergence errors, reduces the number of iterations, and alleviates the communication loads. The SPDS-based decentralized control scheme also enjoys privacy and security features, as well as guaranteed convergence. Additionally, the decentralized control scheme we develop allows chargers to solve their own control problem meaning that (1) device-level constraints, such as charging characteristics that preserve battery state of health, can be met and (2) only a small amount of information needs to be transmitted to a central operator, which protects customers' privacy.

The ultimate goal of this paper is to develop a decentralized EV charging framework that can achieve valley-filling under local EV and distribution network constraints. This framework is built upon a proposed novel decentralized optimization algorithm -- SPDS. The contribution of this paper is three-fold.
\begin{itemize}
\item First, we develop a network-aware EV charging control framework. Specifically, while controlling EVs to achieve valley-filling, heterogeneous charging needs are satisfied and all nodal voltage magnitudes are regulated within the service range. Though we focus on valley-filling and voltage drop, the proposed framework can be flexibly extended to other grid services and include additional network constraints.

\item Second, the proposed control framework only requires a single-layer communication network and requires no additional computing hardware at buses. In addition, customers' privacy can be guaranteed as their SOC information will not be transmitted through the communication network.

\item Third, we develop a new approach to decentralize the EV charging problem that does not require a regularized Lagrangian. To achieve this, we propose an SPDS algorithm that can either be implemented at charging points or utilized by a central operator to perform parallel computing. Execution of SPDS only requires simple Euclidean projections.
\end{itemize}

The rest of this paper is structured as follows. In Section \ref{Problem_formulation}, both of the EV charging model and the distribution network model are constructed to introduce local constraints and network constraints, respectively. Section \ref{Controller_design} first presents the centralized controller, followed by the decentralized control scheme and the SPDS algorithm. Convergency analyses are also provided in Section \ref{Controller_design}. Simulations and result analyses are given in Section \ref{Simulation_results}. Section \ref{Conclusion} concludes this paper and envisions the future work.

\section{EV Charging and Distribution Network} \label{Problem_formulation}

\subsection{EV charging model}
The dynamics of the state-of-charge (SOC) of the $i$th EV can be represented by a first-order discrete-time system as
\begin{equation} \label{raw_charging_dy}
SOC_i(T+1)=SOC_i(T)+\eta_i\Delta t\frac{\bar{P}_i}{\bar{E}_i}u_i(T),
\end{equation}
where $T$ denotes a general discrete time index, $\eta_i$ is the charging efficiency, $\Delta t$ is the sampling time interval, $\bar{P}_i$ is the maximum charging power, and $\bar{E}_i$ is the battery capacity of the $i$th EV. In addition, the charging rate $u_i(T)$ is the control signal which continuously varies in $[0,1]$.

Let $SOC_{i,\text{ini}}$ and $SOC_{i,\text{des}}$ denote the initial SOC and the desired SOC before leaving of the $i$th EV, respectively. Then the total battery energy required by the $i$th EV is
\begin{equation}
E_{i,r}=\bar{E}_i(SOC_{i,\text{des}}-SOC_{i,\text{ini}}).
\end{equation}
At time {$T$}, let the system state $x_i(T)$ denote the energy remaining to be charged to the $i$th EV in order to achieve the total required energy $E_{i,r}$. The charging dynamics can be written as
\begin{equation} \label{E_left_model}
x_i(T+1)= x_i(T)+B_iu_i(T),
\end{equation}
where $B_i=-\eta_i\Delta t \bar{P}_i$. Suppose the $i$th EV is plugged in at time $k_i$, then in order to ensure the driver's charging requirement can be met, the state should satisfy
\begin{equation} \label{heterogeneous_deadline}
x_i(k_i+K_i)=0,
\end{equation} 
where $k_i+K_i$ denotes the designated charging deadline of the $i$th EV. Let $k$ and $k+K$ denote the valley-filling start and end time, respectively. In this paper, we assume that all individually designated charging windows $[k_i,k_i+K_i]$ cover the common charging window $[k,k+K]$. Under this assumption, our control algorithm identifies solutions that deliver all needed energy to vehicles within the common charging deadline, so that there is no need for them to charge outside the common window. The initial SOCs, charging efficiencies, battery capacities, and maximum charging powers are heterogeneous. Let $n$ denote the number of EVs, then an augmented system can be represented as
\begin{equation} \label{state_space_state_all}
x(T+1)=x(T)+\sum_{i=1}^{n}B_{i,c}u_i(T),
\end{equation}
where
\begin{equation} \label{state_space_def}
\begin{aligned}
x(T)&=[x_1(T)~x_2(T)~\cdots~x_n(T)]^{\mathsf{T}} \in \mathbb{R}^{n} , \\
u(T)&=[u_1(T)~u_2(T)~\cdots~u_n(T)]^{\mathsf{T}} \in \mathbb{R}^{n},
\end{aligned}
\end{equation}
and $B_{i,c}\in \mathbb{R}^n$ is the $i$th column of the matrix $B=\bigoplus_{i=1}^nB_i$. $\bigoplus$ and $\oplus$ denote matrix direct sums hereinafter.

Augmenting system \eqref{state_space_state_all} along the valley-filling period $[k,k+K]$, we have
\begin{equation}
\begin{aligned}
\mathcal{X}(k)&=\left[ x(k+1|k)^{\mathsf{T}}~x(k+2|k)^{\mathsf{T}}~\cdots~x(k+K|k)^{\mathsf{T}} \right]^{\mathsf{T}} \\
&= \mathcal{M}x(k)+\sum_{i=1}^{n}\mathcal{B}_{i}\mathcal{U}_i(k)\in\mathbb{R}^{nK},
\end{aligned}
\end{equation}
where
\begin{equation}
\begin{aligned}
\mathcal{M}&=\left[ \bold{I}_n~\bold{I}_n~\cdots~\bold{I}_n\right]^{\mathsf{T}}\in\mathbb{R}^{nK\times n} \\
\mathcal{B}_{i}&=\left[ \begin{array}{cccc}
B_{i,c} &&& \\
B_{i,c} & B_{i,c} & & \\
\vdots & \vdots & \ddots & \\
B_{i,c} & B_{i,c}& \cdots & B_{i,c} \end{array}
\right] \in \mathbb{R}^{nK\times K}, \nonumber
\end{aligned}
\end{equation}
and
\begin{equation}
\mathcal{U}_i(k)=\left[ \begin{array}{c}
u_i(k|k) \\
u_i(k+1|k) \\
\vdots \\
u_i(k+K-1|k)
\end{array}\right] \in \mathbb{R}^{K}. \nonumber
\end{equation}
Herein, $x(k+\kappa|k)$, $\kappa=1,\cdots,K$, is the predicted system state at time $k+\kappa$ based on the known knowledge of the system at time $k$ and all control signals priori to time $k+\kappa$; $u(k+\kappa-1|k)$, $\kappa=1,\cdots,K$, is the predicted control signal at time $k+\kappa-1$ based on the known knowledge of the system at time $k$ and all predicted system states priori to time $k+\kappa$.  

In order to satisfy all drivers' charging requirements by the end of valley-filling, we need to guarantee the $K$th element of $\mathcal{X}(k)$ satisfying
\begin{equation} \label{Charging_requirement_constraint}
\mathcal{X}_K(k)=x(k)+\sum_{i=1}^{n}\mathcal{B}_{i,l}\mathcal{U}_i(k) =\boldsymbol{0},
\end{equation}
where
\begin{equation}
\mathcal{B}_{i,l}=\left[ B_{i,c}~B_{i,c}~\cdots~B_{i,c}\right] \in \mathbb{R}^{n \times K}. \nonumber
\end{equation}
Eq. \eqref{Charging_requirement_constraint} will serve as the local constraints for the controller design.

\subsection{Distribution network model}
In this paper, we consider a radial distribution network, which is a structure commonly adopted in the power system literature and applications. As mentioned in Section \ref{Introduction}, with appropriate assumptions, DistFlow model \cite{Baran_TPD_1989} can be linearized to a linear LinDistFlow model, characterizing a linear relationship between bus power or EV charging power and nodal voltages. In this section, we first describe the LinDistFlow model, then fit it into the valley-filling problem.

\begin{figure}[!htb] \centering
\includegraphics[width=0.46\textwidth]{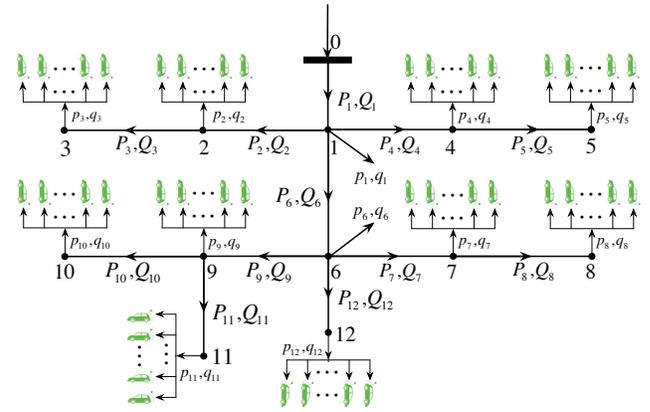}
\caption{A radial distribution network connected with EVs.}
\label{Radial_Network}
\end{figure}

Fig. \ref{Radial_Network} shows the single-phase IEEE-13 Node Test Feeder, in which we discard the transformer between Node 4 and Node 5, discard the switch between Node 6 and Node 7, and assume no capacitor banks. This radial network will be used throughout this paper.
Let $\mathbb{H}=\{\imath| \imath=1,\ldots,h\}$ denote the set of nodes of this distribution feeder and let $\mathbb{E}$ denote the set of all downstream line segments. Node $0$ is the feeder head, decoupling interactions in the downstream distribution system from the rest of the grid and maintaining its own nodal voltage magnitude $\left| V_0 \right|$. To have a clear view of nodal voltage impacts and how controlling EV charging can help alleviate the impacts, we do not consider any possible distributed energy resources, reactive power supplies, or voltage regulations.

At time $T$, let $\left| V_\imath (T) \right|$ denote the voltage magnitude at Node $\imath$; let $p_\imath(T)$ and $q_\imath(T)$ denote the real and reactive power consumption at Node $\imath$; and with a slight abuse of notations, let $r_{\imath \jmath}$ and $x_{\imath \jmath}$ denote the resistance and reactance of the line segment $(\imath, \jmath)$. According to \cite{Baran_TPD_1989, Farivar_CDC_2013, Bansal_CDC_2014}, by omitting the line loss in DistFlow equations, the LinDistFlow model of this distribution network can be written as
\begin{equation} \label{LinDistFlow_original}
\boldsymbol{V}(T)=\boldsymbol{V}_0-2\boldsymbol{R}p(T)-2\boldsymbol{X}q(T),
\end{equation}
where
\begin{equation}
\begin{aligned}
\boldsymbol{V}(T)&=\left[ \begin{array}{c}
\left| V_1(T) \right|^2 \\
\left| V_2(T) \right|^2 \\
\vdots \\
\left| V_h(T) \right|^2
\end{array} \right] \in \mathbb{R}^{h}, \boldsymbol{V}_0=\left[ \begin{array}{c}
\left| V_0 \right|^2 \\
\left| V_0 \right|^2 \\
\vdots \\
\left| V_0 \right|^2
\end{array} \right] \in \mathbb{R}^{h}, \\
p(T)&=\left[ \begin{array}{c}
p_1(T) \\
p_2(T) \\
\vdots \\
p_h(T)
\end{array}\right] \in \mathbb{R}^{h}, q(T)=\left[ \begin{array}{c}
q_1(T) \\
q_2(T) \\
\vdots \\
q_h(T)
\end{array}\right] \in \mathbb{R}^{h}, \nonumber
\end{aligned}
\end{equation}
and
\begin{equation} \label{Definition_R_X}
\begin{aligned}
&\boldsymbol{R}\in\mathbb{R}^{h\times h},~\boldsymbol{R}_{\imath \jmath}=\sum_{(\hat{\imath},\hat{\jmath})\in \mathbb{E}_\imath \cap \mathbb{E}_\jmath}r_{\hat{\imath}\hat{\jmath}}, \\
&\boldsymbol{X}\in\mathbb{R}^{h\times h},~\boldsymbol{X}_{\imath \jmath}=\sum_{(\hat{\imath},\hat{\jmath})\in \mathbb{E}_\imath \cap \mathbb{E}_\jmath}x_{\hat{\imath}\hat{\jmath}}, \nonumber
\end{aligned}
\end{equation}
where $\mathbb{E}_\imath$ and $\mathbb{E}_\jmath$ are the sets containing downstream line segments connecting Node $0$ and Node $\imath$ and connecting Node $0$ and Node $\jmath$, respectively \cite{Farivar_CDC_2013, Bansal_CDC_2014}.

At each node, the power consumption consists of the baseline load (non-EV load) and adjustable load. In this paper, we only consider EV as the adjustable load. Hence, we have
\begin{equation}
\begin{aligned}
p_\imath(T)&=p_{\imath,b}(T)+p_{\imath,EV}(T), \\
q_\imath(T)&=q_{\imath,b}(T)+q_{\imath,EV}(T),
\end{aligned}
\end{equation}
where $p_{\imath,b}(T)$ and $q_{\imath,b}(T)$ denote the real and reactive baseline power at node $\imath$, respectively, and $p_{\imath,EV}$ and $q_{\imath,EV}$ denote the real and reactive EV charging power at node $\imath$, respectively.

Since $-2\boldsymbol{R}p(T)-2\boldsymbol{X}q(T)$ is linear, letting $\boldsymbol{V}_b(T)$ denote the squared voltage drop caused by the baseline load, we have
\begin{equation} \label{V_and_EV}
\boldsymbol{V}(T)=\boldsymbol{V}_0-\boldsymbol{V}_b(T)-2\boldsymbol{R}p_{EV}(T)-2\boldsymbol{X}q_{EV}(T),
\end{equation}
where $p_{EV}(T)=\left[ p_{1,EV}(T)^{\mathsf{T}}~p_{2,EV}(T)^{\mathsf{T}}~\cdots~p_{h,EV}(T)^{\mathsf{T}}\right]^{\mathsf{T}}$ and $q_{EV}(T)=\left[ q_{1,EV}(T)^{\mathsf{T}}~q_{2,EV}(T)^{\mathsf{T}}~\cdots~q_{h,EV}(T)^{\mathsf{T}}\right]^{\mathsf{T}}$.
We further assume EVs only consume real power, resulting in $q_{\imath,EV}(T)=0,~\forall~\imath \in \mathbb{H}$. Thus, \eqref{V_and_EV} is rewritten as
\begin{equation} \label{V_p_EV}
\boldsymbol{V}(T)=\boldsymbol{V}_0-\boldsymbol{V}_b(T)-2\boldsymbol{R}p_{EV}(T).
\end{equation}
Suppose $n_\imath$ EVs are connected at node $\imath$, we have
\begin{equation}
p_{\imath,EV}(T)=\sum_{\hat{\imath}=1}^{n_\imath}\bar{P}_{\imath,\hat{\imath}}u_{\imath,\hat{\imath}}(T),~\imath=1,\ldots,h,
\end{equation}
where $u_{\imath,\hat{\imath}}(T)$ denotes the charging rate of the EV connected at charger $\hat{\imath}$ of node $\imath$, and $\bar{P}_{\imath,\hat{\imath}}$ is the associated maximum charging power. It's worth noting that $\sum_{\imath=1}^{h}n_\imath=n$.

By following ascending orders of $\imath$ and $\hat{\imath}$, we replace the subscripts of $u_{\imath,\hat{\imath}}(T)$ and $\bar{P}_{\imath,\hat{\imath}}$ by $i=1,\ldots,n$. Consequently, we have \eqref{V_p_EV} written as
\begin{equation}
\boldsymbol{V}(T)=\boldsymbol{V}_0-\boldsymbol{V}_b(T)-2\boldsymbol{R}G\bar{P}u(T),
\end{equation}
where
\begin{subequations}
\begin{equation}
G=\bigoplus_{\imath=1}^{h}G_\imath \in \mathbb{R}^{h \times n}, 
\end{equation}
\begin{equation}
\bar{P}=\bigoplus_{i=1}^{n}\bar{P}_i \in \mathbb{R}^{n \times n},
\end{equation}
\end{subequations}
and $G_\imath={\bold{1}}_{n_\imath}^\mathsf{T}$ is the charging power aggregation vector.

Let $D\in \mathbb{R}^{h\times n}$ denote $-2\boldsymbol{R}G\bar{P}$, $y_d(T)$ denote $\boldsymbol{V}_0-\boldsymbol{V}_b(T)$, and $y(T)$ denote $\boldsymbol{V}(T)$, we have
\begin{equation}
y(T)=y_d(T)+Du(T),
\end{equation}
where $y(T)=\left[ y_1(T)~\cdots~y_h(T) \right]^\mathsf{T}\in\mathbb{R}^{h}$, $y_d(T)=\left[ y_{d1}(T)~\cdots~y_{dh}(T) \right]^\mathsf{T}\in\mathbb{R}^h$, and $u(T)$ is defined in \eqref{state_space_def}. Thus, at time $T$, we state the system dynamics as
\begin{subequations}
\begin{equation}
x(T+1)=Ax(T)+Bu(T),
\end{equation}
\begin{equation}
y(T)=y_d(T)+Du(T).
\end{equation}
\end{subequations}

{Augmenting the system output $y(T)$ along the valley-filling period $[k,k+K]$, we have}
\begin{equation}
\mathcal{Y}_k=\mathcal{Y}_{dk}+D_d\left[ \begin{array}{c}
u(k|k) \\
u(k+1|k) \\
\vdots \\
u(k+K-1|k)
\end{array}\right] \in \mathbb{R}^{hK},
\end{equation}
where
\begin{equation}
\begin{aligned}
\mathcal{Y}_k&=\left[ \begin{array}{c}
y(k|k) \\
y(k+1|k) \\
\vdots \\
y(k+K-1|k)
\end{array}\right], \mathcal{Y}_{dk}=\left[ \begin{array}{c}
y_d(k|k) \\
y_d(k+1|k) \\
\vdots \\
y_d(k+K-1|k)
\end{array}\right], \\
D_d&=\underbrace{D\oplus D\oplus \cdots \oplus D}_{K} \in\mathbb{R}^{hK\times nK}. \nonumber
\end{aligned}
\end{equation}
Consequently, it can be obtained that
\begin{equation}
\begin{aligned}
\mathcal{Y}_k&=\mathcal{Y}_{dk}+\left[ \begin{array}{c}
\sum_{i=1}^{n}D_iu_i(k|k) \\
\sum_{i=1}^{n}D_iu_i(k+1|k) \\
\vdots \\
\sum_{i=1}^{n}D_iu_i(k+K-1|k)
\end{array}\right] \\
&=\mathcal{Y}_{dk}+\sum_{i=1}^{n}\mathcal{D}_i\mathcal{U}_i(k),
\end{aligned}
\end{equation}
where
\begin{equation}
\begin{aligned}
\mathcal{D}_i&=\underbrace{D_i \oplus D_i \oplus \cdots \oplus D_i}_{K} \in \mathbb{R}^{hK\times K}, \\
D&=\left[ D_1~D_2~\cdots~D_n\right]. \nonumber
\end{aligned}
\end{equation}

\section{Controller Design} \label{Controller_design}
In the valley-filling problem considered in this paper, the controller is obligated to 1) steer the aggregated EV charging power to fill the overnight valley, 2) guarantee all EVs being charged to their desired SOCs, and 3) maintain the nodal voltage profiles within the service range. In the rest of this section, we first introduce the centralized control problem, then propose a decentralized scheme together with a novel decentralized optimization algorithm. 
\subsection{Centralized controller}
Let $P_b \in \mathbb{R}^K$ denote the aggregated value of distributed uncontrollable loads at all nodes during the whole control period. This $P_b$ is treated as the baseline load, written as
\begin{equation}
P_b=\left[ \begin{array}{c}
P_b(k) \\
P_b(k+1) \\
\vdots \\
P_b(k+K-1)
\end{array}\right] \in \mathbb{R}^K.
\end{equation}
It is assumed that the non-adjustable loads can be well estimated. Designs considering estimation errors will be addressed in our future work. In the rest of this paper, without the loss of generality, we drop the time stamp $k$ for simplicity.

The process of valley-filling is the process of minimizing the variance of the aggregated total load at the feeder, i.e., flattening the total load profile \cite{Gan_TPS_2013}. This implies that the control objective function can be written as
\begin{equation} \label{Centralized_objective}
\begin{aligned}
\mathcal{F}(\mathcal{U})&=\frac{1}{2}\left\| P_b+\sum_{i=1}^{n}f_i(\mathcal{U}_i)\right\|_2^2+\frac{\rho}{2}\left\| \mathcal{U}\right\|_2^2\\
&=\frac{1}{2}\left\| P_b+\tilde{P}\mathcal{U} \right\|_2^2+\frac{\rho}{2}\left\| \mathcal{U}\right\|_2^2,
\end{aligned}
\end{equation}
where $\mathcal{U}=\left[\mathcal{U}_1^\mathsf{T}~\mathcal{U}_2^\mathsf{T}~\cdots~\mathcal{U}_n^\mathsf{T} \right]^\mathsf{T}$, $f_i(\mathcal{U}_i)=\bar{P}_i\mathcal{U}_i$, and $\tilde{P} \in \mathbb{R}^{K\times nK}$ is defined as
\begin{equation}
\tilde{P}=\left[  \underbrace{\bar{P}_1\oplus \bar{P}_1 \oplus \cdots \oplus \bar{P}_1}_{K}~\cdots ~ \underbrace{\bar{P}_n\oplus \bar{P}_n \oplus \cdots \oplus \bar{P}_n}_{K}  \right] . \nonumber
\end{equation}
Note that, the term $\frac{\rho}{2}\left\| \mathcal{U}\right\|_2^2$ in {\eqref{Centralized_objective}} is a proxy of battery degradation cost which can be alternatively represented by a second order polynomial of the charging rates \cite{Ma_TCST_2015}. In other words, given a fixed energy requirement, charging profiles should be as smooth as possible to reduce battery degradation cost.

In this optimal EV charging control problem, there are two types of constraints: local constraints and network constraints. For the $i$th EV, its charging rate should satisfy 
\begin{equation}
\mathcal{U}_i\in \mathbb{U}_i,
\end{equation}
where
\begin{equation}
\mathbb{U}_i := \left\{ \mathcal{U}_i | {\bold{0}}\leq \mathcal{U}_i \leq {\bold{1}}, x_i(k)+\mathcal{B}_{i,l}\mathcal{U}_i=0\right\}. \nonumber
\end{equation}
Note that the set $\mathbb{U}_i$ is convex and it guarantees that, during the controlled charging period, the $i$th EV will reach the desired SOC by varying its charging power in $[0,\bar{P}_i]$.

The cooperative network constraints consider nodal voltage profiles. In order to maintain the power quality, nodal voltage magnitude should be kept within the service range $[\underline{\nu}|V_0|, \overline{\nu}|V_0|]$. Since no DER, reactive power supply, or voltage regulation is considered, only a lower bound as follows is necessary for bus voltage magnitudes.
\begin{equation} \label{lower_bound_voltage}
\begin{aligned}
\mathcal{Y}_k&=\mathcal{Y}_{dk}+\sum_{i=1}^{n}\mathcal{D}_i\mathcal{U}_i \\
&\geq {\underline{\nu}}^2 \boldsymbol{V}_0.
\end{aligned}
\end{equation}
Herein, we assume that with the baseline load bus voltage magnitudes are within $[\underline{\nu}|V_0|, \overline{\nu}|V_0|]$ and that voltages will be less than or equal to the baseline during EV charging. This allows us to only consider lower voltage constraints in EV charging decisions. Let $\mathcal{Y}_b$ denote ${\underline{\nu}}^2 \boldsymbol{V}_0-\mathcal{Y}_{dk}$, the cooperative network constraint can be represented as
\begin{equation}
\mathcal{Y}_b-\sum_{i=1}^{n}\mathcal{D}_i\mathcal{U}_i \leq {\bold{0}}.
\end{equation}

To summarize, the optimal control sequences for all connected EVs providing valley-filling can be obtained by solving
\begin{equation} \label{Centralized_Problem}
\begin{aligned}
\min_{\mathcal{U}} &~ \mathcal{F}(\mathcal{U}) \\
\text{s.t.}&~\mathcal{U}_i \in \mathbb{U}_i,~\forall~i=1,2,\ldots,n, \\
&~\mathcal{Y}_b-\sum_{i=1}^{n}\mathcal{D}_i\mathcal{U}_i \leq {\bold{0}}.
\end{aligned}
\end{equation}
in a centralized fashion.

\subsection{Decentralized control scheme}
It can be clearly observed from \eqref{Centralized_objective} that, individual $\mathcal{U}_i$'s in the objective function $\mathcal{F}(\mathcal{U})$ are coupled and non-separable by the 2-norm and the constraints in \eqref{Centralized_Problem} contain a linearly coupled inequality constraint. This can be understood as that all EVs cooperate to achieve valley-filling, at the meantime their charging powers are jointly constrained by the nodal voltages. Distributing the computational load to individual EV chargers and realizing a decentralized control scheme for this problem is challenging.

Let
\begin{equation}
\mathcal{F}\left( \mathcal{U}\right) =\mathcal{G}\left( \mathcal{U}\right)+\frac{\rho}{2}\left\| \mathcal{U}\right\|_2^2.
\end{equation}
Then, the Lagrange dual problem of \eqref{Centralized_Problem} is
\begin{equation} \label{Dual_Problem}
\max_{\lambda\in \mathbb{R}_{+}^{hK}}\left\{ \min_{\mathcal{U}\in \mathbb{U}} \mathcal{L}(\mathcal{U},\lambda)\right\},
\end{equation}
where $\mathbb{U}=\mathbb{U}_1\times \mathbb{U}_2 \times \cdots \times \mathbb{U}_n$, $\lambda \in \mathbb{R}^{hK}_+$ is the dual variable associated with the inequality constraint $\mathcal{Y}_b-\sum_{i=1}^{n}\mathcal{D}_i\mathcal{U}_i \leq {\bold{0}}$, and $\mathcal{L}:\mathbb{R}^{nK}\times\mathbb{R}^{hK}\rightarrow\mathbb{R}^1$ is the Lagrangian given by
\begin{equation} \label{Centralized_Lag}
\begin{aligned}
\mathcal{L}(\mathcal{U},\lambda)&=\mathcal{G}(\mathcal{U})+\frac{\rho}{2}\left\| \mathcal{U}\right\|_2^2+\lambda^\mathsf{T}\left( \mathcal{Y}_b-\sum_{i=1}^{n}\mathcal{D}_i\mathcal{U}_i \right) \\
&=\mathcal{G}(\mathcal{U})+\frac{\rho}{2}\left\| \mathcal{U}\right\|_2^2+\lambda^\mathsf{T}d(\mathcal{U}).
\end{aligned}
\end{equation}

Before touching the decentralized algorithm, we first introduce a useful theorem.

{\bf{Theorem 1:}} \emph{Saddle-Point Theorem \cite{Boyd_book_2004}: The point $(\mathcal{U}^*,\lambda^*)\in\mathbb{U}\times\mathbb{R}^{hK}_+$ is a primal-dual solution of problems \eqref{Centralized_Problem} and \eqref{Dual_Problem} if and only if there holds}
\begin{equation}
\mathcal{L}(\mathcal{U}^*,\lambda) \leq \mathcal{L}(\mathcal{U}^*,\lambda^*) \leq \mathcal{L}(\mathcal{U},\lambda^*),~ \forall~(\mathcal{U},\lambda) \in \mathbb{U} \times \mathbb{R}_+^{hK}. ~~~~~~\blacksquare \nonumber
\end{equation} 

{\bf{Theorem 1}} indicates that, given that the Lagrangian $\mathcal{L}(\mathcal{U},\lambda)$ is a convex-concave function defined on $\mathbb{U}\times \mathbb{R}_{+}^{hK}$, the solutions to the saddle-point problem \eqref{Dual_Problem} solve the original problem \eqref{Centralized_Problem}. Further, the optimal solution $(\mathcal{U}^*,\lambda^*)$ to \eqref{Dual_Problem} can be obtained by solving a system of projection equations \cite{Koshal_SIAM_2011}
\begin{subequations} \label{Centralized_Projection}
\begin{equation}
\mathcal{U}^*=\Pi_{\mathbb{U}}\left( \mathcal{U}^*-\nabla_{\mathcal{U}}\mathcal{L}(\mathcal{U}^*,\lambda^*)\right),
\end{equation}
\begin{equation}
\lambda^*=\Pi_{\mathbb{R}^{hK}_+}\left( \lambda^*+\nabla_{\lambda}\mathcal{L}(\mathcal{U}^*,\lambda^*)\right),
\end{equation}
\end{subequations}
where $\Pi_\mathbb{X}({\bold{x}})$ is the projection function projecting $\bold{x}$ into the convex set $\mathbb{X}$. Thus, $\mathcal{U}^*$ is an optimal solution of \eqref{Centralized_Problem} if and only if it is a solution to the system \eqref{Centralized_Projection} for some $\lambda^* \in \mathbb{R}^{hK}_+$.

The classic primal-dual subgradient (PDS) algorithm \cite{Arrow_book_1958} is a good candidate for solving \eqref{Centralized_Projection} in a centralized way and it suggests the $\ell$th iteration of updating be performed as
\begin{subequations}
\begin{equation}
\mathcal{U}^{(\ell+1)}=\Pi_{\mathbb{U}}\left( \mathcal{U}^{(\ell)}-\alpha_\ell\nabla_{\mathcal{U}}\mathcal{L}(\mathcal{U}^{(\ell)},\lambda^{(\ell)})\right),
\end{equation}
\begin{equation}
\lambda^{(\ell+1)}=\Pi_{\mathbb{R}^{hK}_+}\left( \lambda^{(\ell)}+\beta_\ell\nabla_{\lambda}\mathcal{L}(\mathcal{U}^{(\ell)},\lambda^{(\ell)})\right),
\end{equation}
\end{subequations}
where $\alpha_\ell>0$ and $\beta_\ell>0$ are the iteration step sizes for the primal and dual variables, respectively. 

To consider the problem in a decentralized manner, let $\mathcal{U}_{-i}$ denote the collection of all $\mathcal{U}_j$, $j\neq i$. Assume the problem \eqref{Centralized_Problem} satisfies the Slater condition, i.e., 
\begin{equation}
\exists~\bar{\mathcal{U}}\in \mathbb{U},~\text{such that}~\mathcal{Y}_b-\sum_{i=1}^{n}\mathcal{D}_i\bar{\mathcal{U}}_i \leq {\bold{0}},
\end{equation}
where $\bar{\mathcal{U}}$ is the Slater point, the strong duality holds. Define the variational inequality as: Given a subset $\mathbb{K}\subseteq \mathbb{R}^n$ and a mapping $\mathcal{P}:\mathbb{K}\mapsto \mathbb{R}^n$, the variational inequality, denoted by $\text{VI}(\mathbb{K},\mathcal{P})$, is to find $\bold{x}\in\mathbb{K}$ such that
\begin{equation}
(\bold{y}-\bold{x})^{\mathsf{T}}\mathcal{P}(\bold{x})\geq 0,~\forall~\bold{y}\in\mathbb{K}.
\end{equation}
Then, from the first-order optimality conditions and the decomposable structure of $\mathbb{U}$, it can be seen that $(\mathcal{U}^*,\lambda^*)\in\mathbb{U}\times \mathbb{R}_+^{hK}$ is a solution to \eqref{Centralized_Problem} if and only if $\mathcal{U}_i^*$ solves the parameterized variational inequality $\text{VI}\left(\mathbb{U}_i,\nabla_{\mathcal{U}_i}\mathcal{L}(\mathcal{U}_i;\mathcal{U}_{-i}^*,\lambda^*) \right)$, $i=1,\ldots,n$, and $\lambda^*$ solves $\text{VI}\left( \mathbb{R}_+^{hK},-\nabla_{\lambda}\mathcal{L}(\mathcal{U}^*,\lambda)\right)$. In addition, $(\mathcal{U}^*,\lambda^*)$ solves $\text{VI}\left(\mathbb{U}_i,\nabla_{\mathcal{U}_i}\mathcal{L}(\mathcal{U}_i;\mathcal{U}_{-i}^*,\lambda^*) \right)$, $i=1,\ldots,n$, and $\text{VI}\left( \mathbb{R}_+^{hK},-\nabla_{\lambda}\mathcal{L}(\mathcal{U}^*,\lambda)\right)$ if and only if each $\mathcal{U}_i^*$ is a zero of the parameterized natural map ${\bold{F}}_{\mathbb{U}_i}^{\text{nat}}\left( \mathcal{U}_i;\mathcal{U}_{-i}^*,\lambda^* \right)=0$, $i=1,\ldots,n$, and $\lambda^*$ is a zero of the parameterized natural map ${\bold{F}}_{\mathbb{R}_+^{hK}}^{\text{nat}}\left( \lambda; \mathcal{U}^*\right)=0$ \cite{Koshal_SIAM_2011,Facchinei_book_2003}, where the parametrized natural maps are defined as
\begin{equation}
\begin{aligned}
{\bold{F}}_{\mathbb{U}_i}^{\text{nat}}\left( \mathcal{U}_i;\mathcal{U}_{-i}^*,\lambda^* \right)&\triangleq  \mathcal{U}_i-\Pi_{\mathbb{U}_i}\left(\mathcal{U}_i-\nabla_{\mathcal{U}_i}\mathcal{L}\left(\mathcal{U}_i;\mathcal{U}_{-i}^*,\lambda^*\right)\right),\\
{\bold{F}}_{\mathbb{R}_+^{hK}}^{\text{nat}}\left( \lambda; \mathcal{U}^*\right)& \triangleq  \lambda -\Pi_{\mathbb{R}_+^{hK}}\left(\lambda+\nabla_{\lambda}\mathcal{L}\left( \mathcal{U}^*,\lambda \right) \right). \nonumber
\end{aligned}
\end{equation}

Based on the parametrized natural maps, Koshal \emph{et al.} \cite{Koshal_SIAM_2011} developed a strategy to decentralize the problem via a RPDS algorithm. The regularization terms in the Lagrangian ensure convergence, however they also introduce errors relative to the optimal solution of the original problem. In this paper we explore an alternative approach -- which we call a shrunken-primal-dual subgradient (SPDS) algorithm -- to circumvent the need for a regularization term in the Lagrangian while still solving the problem in a decentralized manner. The SPDS works as follows: At the $\ell$th iteration, primal and dual variables update by following
\begin{subequations} \label{all_updates}
\begin{equation} \label{primal_updates}
\mathcal{U}_i^{(\ell+1)} =\Pi_{\mathbb{U}_i}\left(\frac{1}{\tau_\mathcal{U}}\Pi_{\mathbb{U}_i}\left( \tau_\mathcal{U}\mathcal{U}_i^{(\ell)}-\alpha_{(i,\ell)}\nabla_{\mathcal{U}_i}\mathcal{L}(\mathcal{U}^{(\ell)},\lambda^{(\ell)})\right)\right),
\end{equation}
\begin{equation} \label{dual_updates}
\lambda^{(\ell+1)}=\Pi_{\mathbb{D}}\left(\frac{1}{\tau_\lambda}\Pi_{\mathbb{D}}\left( \tau_\lambda\lambda^{(\ell)}+\beta_\ell\nabla_{\lambda}\mathcal{L}(\mathcal{U}^{(\ell)},\lambda^{(\ell)})\right)\right),
\end{equation}
\end{subequations}
where
\begin{equation}
\mathbb{D}:=\left\{ \lambda| \lambda\geq {\bold{0}},~\|\lambda\|_2 \leq d_\lambda  \right\},
\end{equation}
and $0<\tau_\mathcal{U},\tau_\lambda<1$. 

The proposed SPDS features a two-tier projection in which primal and dual variables are shrunken and then expanded. Take \eqref{primal_updates} for example, at tier-1 projection, solution from the previous iteration $\mathcal{U}_i^{(\ell)}$ is shrunken by $\tau_\mathcal{U}$, moved towards the descent direction of the Lagrangian, and then projected into $\mathbb{U}_i$. Because the shrinkage brings a certain level of conservativeness, result from tier-1 projection is first expanded by $1/\tau_\mathcal{U}$ and then projected back into $\mathbb{U}_i$ at the tier-2 projection. The shrinkage-expansion and the two-tier projection in \eqref{all_updates} guarantee the convergency without any regularization. In the proposed SPDS,
\begin{equation}
d_\lambda \triangleq \frac{\mathcal{F}(\bar{\mathcal{U}})-\tilde{l}}{\gamma}+\sigma,
\end{equation}
where $\bar{\mathcal{U}}$ is a Slater point of \eqref{Centralized_Problem}, $\sigma>0$, $\gamma=\min_{j=1,\ldots,hK}\left\{-d_j(\bar{\mathcal{U}})\right\}$, and
\begin{equation}
\tilde{l}=\min_{\mathcal{U}_i\in \mathbb{U}_i,i=1,\ldots,n}\mathcal{L}(\mathcal{U}_i, \ldots,\mathcal{U}_n,\tilde{\lambda}),~\forall\tilde{\lambda}\in \mathbb{R}^{hK}_+,
\end{equation}
where $d_j(\mathcal{U})$ is the $j$th entry of $d(\mathcal{U})$. It is shown in \cite{Chang_TAC_2014, Nedic_JOTA_2009, Nedic_SJO_2009} that, under Slater condition, the primal solution $\mathcal{U}_i^*$ exists and therefore strong duality holds. Further, the dual optimal set $\mathbb{D}$ is not empty and bounded by $d_\lambda$. In the rest of this paper, we assume that the bound $d_\lambda$ can be obtained \emph{a priori}.

The proposed SPDS algorithm is summarized in {\bf{Algorithm \ref{Decentralized_Alg}}}.
\begin{algorithm}[!htb]
\caption{SPDS Algorithm}\label{Decentralized_Alg}
\begin{algorithmic}[1]
\State Iteration number $\ell=0$; EVs initialize $\mathcal{U}_i^{(0)}$; Operator initializes $\lambda^{(0)}$; Tolerance $\tau_\epsilon$; Initial error $\epsilon=10^9$; Maximum iteration $\ell_{max}$; 
\Procedure {~}{}
\While {$\epsilon>\tau_\epsilon$ and $\ell\leq \ell_{max}$}
\State Each EV charger transmits its own $\mathcal{U}_i^{(\ell)}$ to the operator.
\State Operator computes the Lagrangian gradient and broadcasts it together with $\lambda^{(\ell)}$ to all chargers. 
\State All connected chargers perform \eqref{primal_updates}.
\State The operator performs \eqref{dual_updates}.
\State $\epsilon=\left\| \mathcal{U}^{(\ell+1)}-\mathcal{U}^{(\ell)} \right\|_2$.
\State $\ell=\ell+1$.
\EndWhile
\EndProcedure 
\end{algorithmic}
\end{algorithm}
The merit of this algorithm is that no communication network is needed between the chargers and all chargers can perform the computing in a parallel fashion.

{\bf{Remark 1:}} In {\bf{Algorithm 1}}, individual EV chargers only need to transmit their calculated control sequences to the operator for updating the dual variable. Hence, drivers' private information, especially the SOC, can be well protected, complying with both ISO/IEC 15118 and SAE 2847.
\hfill $\blacksquare$ 

\subsection{SPDS Convergency analysis}

In this paper, we will only be discussing the constant step size, i.e., $\alpha_{i,\ell}=\alpha$ and $\beta_{\ell}=\beta$. Adaptive step sizes will be studied in our future work. Having constant step sizes, the updates of primal and dual variables in SPDS are represented as
\begin{subequations} \label{constant_update}
\begin{equation}
\mathcal{U}_i^{(\ell+1)} =\Pi_{\mathbb{U}_i}\left(\frac{1}{\tau_\mathcal{U}}\Pi_{\mathbb{U}_i}\left( \tau_\mathcal{U}\mathcal{U}_i^{(\ell)}-\alpha\nabla_{\mathcal{U}_i}\mathcal{L}(\mathcal{U}^{(\ell)},\lambda^{(\ell)})\right)\right),
\end{equation}
\begin{equation}
\lambda^{(\ell+1)}=\Pi_{\mathbb{D}}\left(\frac{1}{\tau_\lambda}\Pi_{\mathbb{D}}\left(\tau_\lambda \lambda^{(\ell)}+\beta\nabla_{\lambda}\mathcal{L}(\mathcal{U}^{(\ell)},\lambda^{(\ell)})\right)\right).
\end{equation}
\end{subequations}

We then present the following convergency theorem.

{\bf{Theorem 2:}} \emph{Let $\{\zeta^{(\ell)}\}$, where $\zeta^{(\ell)}=\left[ {\mathcal{U}^{(\ell)}}^{\mathsf{T}}~{\lambda^{(\ell)}}^{\mathsf{T}}\right]^{\mathsf{T}}$, be a sequence generated by \eqref{constant_update}. Then we have
\begin{equation} \label{zeta_convergency}
\left\| \zeta^{(\ell+1)}-\zeta^* \right\|_2 \leq \sqrt{\varrho} \left\| \zeta^{(\ell)}-\zeta^* \right\|_2,~\forall~\ell\geq 0,
\end{equation}
where
\begin{equation}
\varrho=\max\left\{ \frac{\hat{\alpha}}{\alpha},\frac{\hat{\beta}}{\beta}\right\}+\delta\max\{\hat{\alpha},\hat{\beta}\}L_{\Phi}^2-2c\min\{\hat{\alpha},\hat{\beta}\}+|\hat{\tau}|\phi. \nonumber
\end{equation}
Herein, $\hat{\alpha}=\alpha/\tau_\mathcal{U}^2$, $\hat{\beta}=\beta/\tau_\lambda^2$, $\hat{\tau}=\hat{\alpha}-\hat{\beta}$, and
\begin{equation}
\delta=\left\{ \begin{array}{cc}
\alpha & \text{if}~\hat{\alpha}>\hat{\beta},\\ 
\beta & \text{if}~\hat{\alpha}<\hat{\beta}.
\end{array}\right. \nonumber
\end{equation}
In addition,
\begin{equation}
\begin{aligned}
\phi=&\max \left\{ L_d^2,1-(1-\sgn(\hat{\tau}))\frac{1-\tau_\lambda}{\beta} \right. \\
&~~~~~~~~~~~~~\left. -(1+\sgn(\hat{\tau}))\left( \rho+\frac{1-\tau_\mathcal{U}}{\alpha} \right) \right\}.
\end{aligned}
\end{equation}
$L_{\Phi}$ is the Lipschitz constant for the mapping
\begin{equation}
\Phi(\zeta)= \left[ \begin{array}{c}
\nabla_\mathcal{U}\mathcal{L}(\mathcal{U},\lambda)+\frac{1-\tau_\mathcal{U}}{\alpha}\mathcal{U} \\
-\nabla_\mathcal{\lambda}\mathcal{L}(\mathcal{U},\lambda)+\frac{1-\tau_\lambda}{\beta}\lambda
\end{array}\right]=\left[ \begin{array}{c}
\Phi_1(\zeta) \\
\Phi_2(\zeta)
\end{array}\right],
\end{equation}
and $c$ is the strongly monotone constant of the mapping $\Phi(\zeta)$. Furthermore, if the tuple $(\alpha,\beta)$ satisfies
\begin{equation} \label{theorem_2_con}
\max\left\{ \frac{1}{\tau_\mathcal{U}^2}, \frac{1}{\tau_\mathcal{\lambda}^2}\right\}+\delta\max\left\{ \hat{\alpha},\hat{\beta}\right\}L_{\Phi}^2-1<2c\max\left\{ \hat{\alpha},\hat{\beta}\right\},
\end{equation}
$\varrho$ is guaranteed within $(0,1)$, indicating that the sequence $\{\zeta^{(\ell)}\}$ converges as $\ell \rightarrow \infty$.} \hfill $\blacksquare$

\emph{Proof of {\bf{Theorem 2}}:} Please see the \textsc{Appendix}. 

{\bf{Remark 2:}} Both the tuples $(\alpha,\beta)$ and $(\tau_{\mathcal{U}}, \tau_{\lambda})$ affect the SPDS convergence and its speed. Given a tuple $(\tau_{\mathcal{U}},\tau_\lambda)$, the tuple $(\alpha,\beta)$ can either be theoretically chosen by satisfying \eqref{theorem_2_con} or be empirically tuned up from small values to guarantee and accelerate the convergence. Shrinking parameters $\tau_{\mathcal{U}}$ and $\tau_{\mathcal {\lambda}}$ should be selected by following the principles below. Firstly, the fundamental principle is $\tau_{\mathcal{U}}, \tau_{\mathcal {\lambda}} \in (0,1)$. If this is violated, e.g., $\tau_{\lambda}=1$, the strongly monotone constant $c$ as defined in \eqref{monotone_constant} becomes 0, immediately violating the sufficient convergence condition \eqref{theorem_2_con} and resulting in an impossible selection of tuple $(\alpha,\beta)$. Secondly, $\tau_{\mathcal{U}}$ and $\tau_{\mathcal {\lambda}}$ should not be too close to 1. Otherwise, value of the strongly monotone constant $c$ gets close to 0, which can also lead to the violation of \eqref{theorem_2_con}. Thirdly, small $\tau_{\mathcal{U}}$ and $\tau_{\mathcal {\lambda}}$ can always guarantee the convergence. However, the convergence speed in terms of the value of $1/\varrho$ is monotonically decreasing as $\omega=\max\left\{ \frac{1}{\tau_\mathcal{U}^2}, \frac{1}{\tau_\mathcal{\lambda}^2}\right\}$ increases. Thus, we always want to choose $\tau_{\mathcal{U}}$ and $\tau_{\mathcal {\lambda}}$ as large as possible to accelerate the convergence. Hence, $\tau_{\mathcal{U}}$ and $\tau_\lambda$ should be chosen close to 1 but not too close. In the later simulations, both $(\alpha,\beta)$ and $(\tau_{\mathcal{U}},\tau_{\lambda})$ are empirically chosen and tuned. We are currently investigating on an analytical approach of selecting $(\tau_{\mathcal{U}}, \tau_\lambda)$ and will present it in our future work. \hfill $\blacksquare$


\section{Simulation Results} \label{Simulation_results}
Fig. \ref{Radial_Network} depicts the residential distribution network used in our simulations. Note that, Node $1$ and Node $6$ have no EV connected, and each of the other nodes is connected with $70$ houses equipped with level-2 chargers (maximum 6.6 kW). Battery capacities are uniformly distributed in $[18, 20]$ kWh. Initial and designated SOCs are uniformly distributed in $[0.3,0.5]$ and $[0.7,0.9]$, respectively. Primal and dual step sizes are empirically tuned to $\alpha_{(i,\ell)}=\alpha=2.8\times10^{-10}$ and $\beta_{\ell}=\beta=1.8$, respectively; $d_\lambda$ is selected as $d_\gamma=5\times 10^5$; shrinking parameters are empirically chosen as $\tau_\mathcal{U}=\tau_\lambda=0.974$. Maximum iteration number $\ell_{max}=25$; tolerance $\tau_\epsilon=1\times 10^{-4}$. The above parameters were carefully chosen to accelerate the convergence upon complying with the convergence conditions. Note that the parameters should be tuned if a different network model is given. Initial values of $\mathcal{U}_i^{(0)}(k)$ and $\lambda^{(0)}$ are both set to all-zero vectors. We set the voltage lower bound to $\underline{\nu}=0.954$, which is a bit higher than the ANSI C84.1 service standard, to accommodate the discarded line losses in LinDistFlow model. The baseline load data is collected from Southern California Edison \cite{SCE_reg} and scaled to fit the simulations. The valley-filling service period is from 19:00 to 8:00 next day. All simulations are conducted in MATLAB+cvx on a MacBook Pro with 2.8 GHz Intel Core i7 and 16 GB memory.

\subsection{Centralized control}
\begin{figure}[!htb] \centering
\includegraphics[width=0.47\textwidth, trim = 3mm 0mm 15mm 0mm, clip]{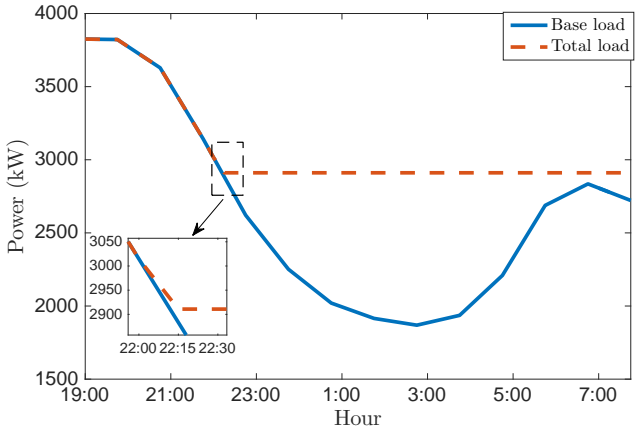}
\caption{Baseline load and total load under centralized control.}
\label{P_agg_cen}
\end{figure}

\begin{figure}[!htb] \centering
\includegraphics[width=0.47\textwidth, trim = 8mm 3mm 18mm 0mm, clip]{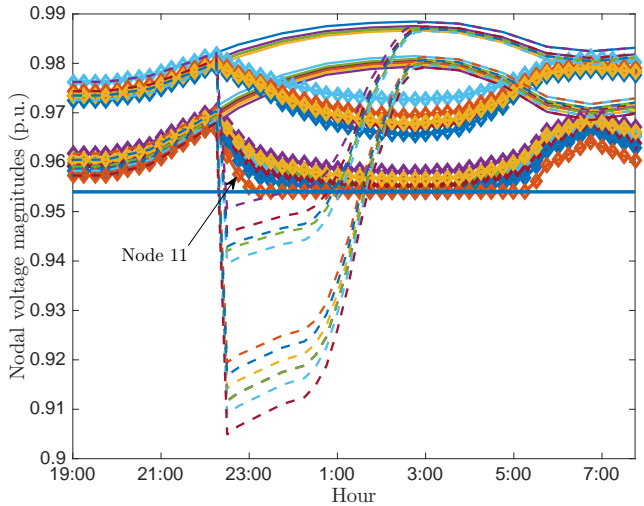}
\caption{Nodal voltage magnitudes, computing via the LinDistFlow equations, of baseline load (solid lines), total load under centralized control (diamond-marked lines), and uncontrolled total load (dashed lines).  }
\label{Voltage_cen}
\end{figure}

By applying the control signals obtained from solving \eqref{Centralized_Problem} in the centralized way, we can find the comparisons between the aggregated total load and baseline load in Fig. \ref{P_agg_cen}.
 Specifically, all EV chargers are commanded to be idle before 22:00; EVs start charging at 22:00 to fill the valley. It can be observed that the aggregated total load profile is not strictly flat after EVs start charging; however its value decays smoothly to a flat one during 22:00--22:30. This transition is a result of the nodal voltage constraints and local charging constraints. Though the profile is not strictly flat, it is the optimal performance. To have a better understanding of this, we show the nodal voltage magnitudes of the baseline load and total loads in Fig. \ref{Voltage_cen},
in which solid lines, diamond-marked lines, and dashed lines represent the nodal voltage magnitudes of baseline loads, controlled total loads, and uncontrolled total loads, respectively (voltages are computed via the LinDistFlow equations; we compare results to DistFlow calculations in Section~\ref{sec:model_discrepancy}). Following observations can be made:
\begin{itemize}
\item In the controlled case, the voltage magnitude of Node 11 is in line with the 0.954 p.u. during 23:15--5:15.
\item In the controlled case, voltage magnitudes of Nodes 8 and 10 are slightly above 0.954 p.u. during 1:15--4:00.
\item Without an appropriate control, nodal voltages will drop far below the lower bound.
\end{itemize}
The above observations reveal that, in the controlled case, there is no room for the total load to increase after 23:15, otherwise, the nodal voltage magnitudes cannot be maintained. In addition, since the optimization objective is to minimize the variance of the total load profile, considering the integral of the difference between the total and baseline loads should be a constant and the fact that the total load profile before 22:00 is higher than the flat value, the optimal profile must have the transition from high to flat as early and as smooth as possible. 

\begin{figure}[!htb] \centering
\includegraphics[width=0.47\textwidth, trim = 3mm 5mm 6mm 10mm, clip]{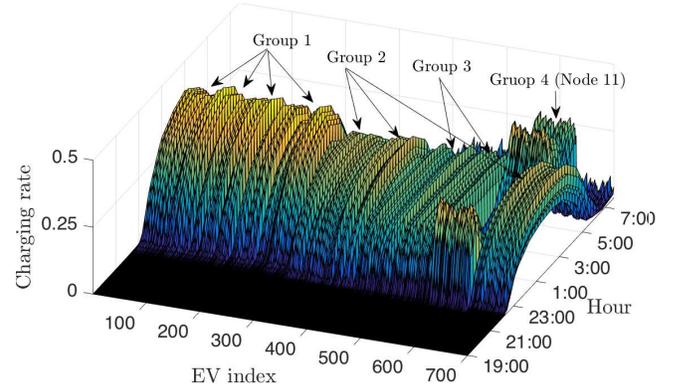}
\caption{Charging profiles of all EVs under centralized control.}
\label{Rate_Con_cen}
\end{figure}

Charging profiles of all EVs under the centralized control are shown in Fig. \ref{Rate_Con_cen}.
Each slice vertical to the EV index axis represents the charging profile of a specific EV during the service period. The 700 EVs are indexed with an ascending order of the node index. It is clear that the charging profiles of EVs are not similar. Classify those EVs into four groups, i.e., Group 1 (Nodes 2 to 5), Group 2 (Nodes 7, 8, and 12), Group 3 (Nodes 9 and 10), and Group 4 (Node 11), according to their geographic locations and downstream line segment impedance, we can notice that charging profiles of EVs in the same group are similar. Geographically, EVs in Group 1 are closer to the feeder head than other EVs, indicating that the network nodal voltages are less sensitive to their charging powers. At around 3:00 when the power need from the EVs for valley-filling reaches the maximum, only EVs having less voltage sensitivities can be pushed to high charging rates. Similarly, since EVs in Group 4 are the furthest to the feeder head and the downstream line impedance are the largest, the network nodal voltages are much more sensitive to their charging powers. Thus, it is optimal to keep the charging rates of EVs in Group 4 as low as possible to compensate for large voltage drops caused by high power needs. EVs in Groups 2 and 3 are in between, thus their charging profiles are rational to be parabola-like or trapezoid-like with lower peak values.


It's worth mentioning that all EVs are charged to the designated SOCs under the centralized control. Figures showing the SOC evolutions are omitted here due to the space limit.

\subsection{Decentralized control}

\begin{figure}[!htb] \centering
\includegraphics[width=0.5\textwidth, trim = 3mm 1mm 15mm 3mm, clip]{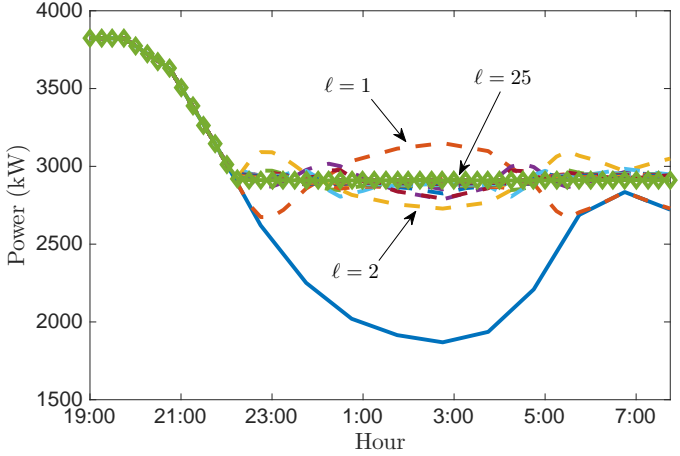}
\caption{Baseline load and total loads in 25 iterations under decentralized control.}
\label{P_agg_Decen}
\end{figure}

Fig. \ref{P_agg_Decen} shows the baseline load together with the total load under the decentralized control in 25 iterations.
Though 25 iterations are performed, the total load converges to an acceptable profile at the 15th iteration. The normed error between total loads at the 15th and 25th iterations is 0.05\%. Slight adjustments after the 15th iteration are to bind and unbind voltage constraints, further improving the valley-filling performance.

\begin{figure}[!htb] \centering
\includegraphics[width=0.5\textwidth, trim = 8mm 4mm 15mm 3mm, clip]{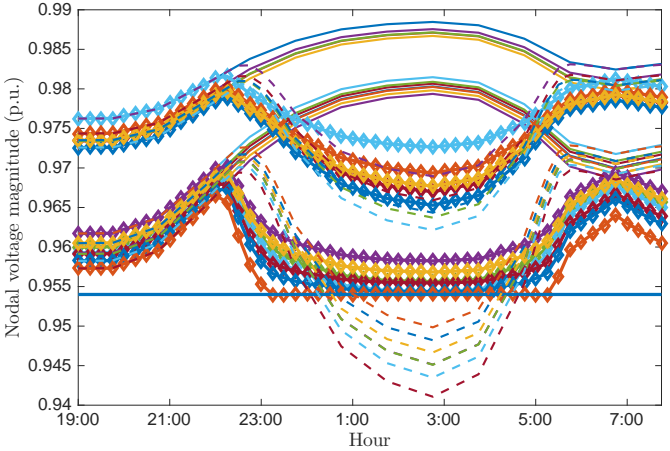}
\caption{Nodal voltage magnitudes of baseline load (solid lines) and total loads at the 1st iteration (dashed lines) and the 25th iteration (diamond-marked lines). Each line represents one node.}
\label{Voltage_Decen}
\end{figure}

Nodal voltage magnitudes of the baseline load, total loads at the 1st iteration and the 25th iteration are shown in Fig. \ref{Voltage_Decen}.
At the 1st iteration (dashed lines), we can clearly observe a large voltage drop below the 0.954 p.u. threshold. This is because at this early iteration number, voltage constraints have yet to be fully addressed and all EVs have parabola-like charging profiles. The voltage converging process is much like a damped oscillation and the voltage profiles finally converge to ones satisfying the constraints.

\begin{figure}[!htb] \centering
\includegraphics[width=0.5\textwidth, trim = 4mm 3mm 4mm 3mm, clip]{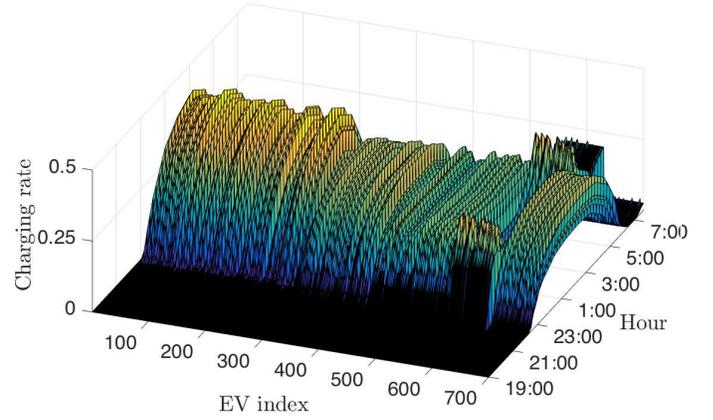}
\caption{Charging profiles of all EVs at the 25th iteration under the decentralized control.}
\label{Rate_Decen}
\end{figure}

\begin{figure}[!htb] \centering
\includegraphics[width=0.5\textwidth, trim = 1mm 5mm 5mm 10mm, clip]{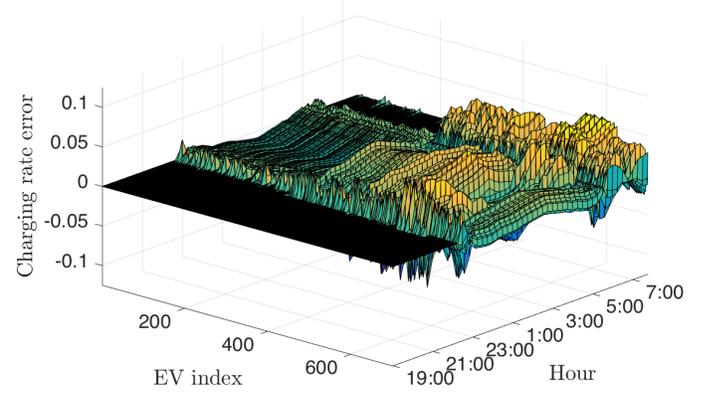}
\caption{The difference between the charging profiles obtained from centralized and decentralized controls..}
\label{rate_diff}
\end{figure}

Charging profiles of all EVs at the 25th iteration are shown in Fig. \ref{Rate_Decen}.
This profile essentially coincides with the one obtained from the centralized controller, but not exactly the same. Errors between them (centralized minus decentralized) are shown in Fig. \ref{rate_diff}.
In most cases, errors are of small magnitudes; large errors (up to 0.05) mainly exist on EVs in Group 4 which are the deciding factors of voltage constraints. It is worth mentioning that, without the stopping criteria, the algorithm can continue adjusting charging rates to converge to the optimal solution. 

\begin{figure}[!htb] \centering
\includegraphics[width=0.5\textwidth, trim = 3mm 2mm 12mm 3mm, clip]{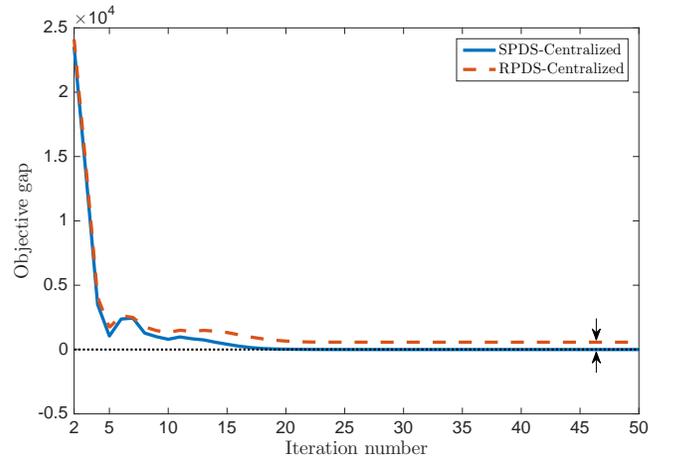}
\caption{Objective value gaps between SPDS and RPDS.}
\label{SPDS_RPDS}
\end{figure}

As aforementioned in Section \ref{Introduction}, one of the design objectives of the SPDS is to eliminate the regularization error introduced by the dual regularization term in RPDS \cite{Koshal_SIAM_2011}. To verify the error-free feature of SPDS, we conducted simulations using RPDS to contrast. Fig. \ref{SPDS_RPDS} shows the objective gaps at different iterations generated by SPDS and RPDS.
Primal and dual step sizes of RPDS were chosen the same as those of SPDS; the coefficient of the dual regularization term in RPDS was chosen as 0.1 to guarantee convergence. It can be readily observed from Fig. \ref{SPDS_RPDS} that SPDS has relatively faster convergence speed and no objective gap exists once converged. In contrast, RPDS converges slightly slower with a certain objective gap, which is caused by the dual regularization term in the Lagrangian. Convergence of primal variables can also be revealed from Fig. \ref{SPDS_RPDS}.

\begin{figure}[!htb] \centering
\includegraphics[width=0.5\textwidth, trim = 3mm 2mm 12mm 3mm, clip]{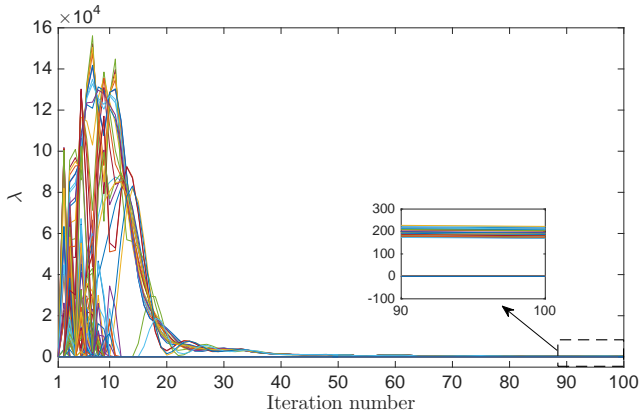}
\caption{Convergency of the dual variable $\lambda$.}
\label{lambda_con_Decen}
\end{figure}

The evolution of the dual variable $\lambda$ in 100 iterations is shown in Fig. \ref{lambda_con_Decen}.
During this process, voltage constraints are being bound and unbound to reach the optimal valley-filling performance. At the final stage, all $\lambda$ values that are above 0 correspond to Node 11, indicating only the voltage constraints at Node 11 are binding. It is known that the sensitivity of all nodal voltage magnitudes to the charging power at node $\jmath$ can be reflected by the $\jmath$th column of the matrix $\boldsymbol{R}$ and is directly determined by the downstream line impedance. Since Node 11 possesses the largest downstream impedance and its voltage is also largely affected by the power at Node 10, the voltage constraint on Node 11 is more likely to be binding. 

It's worth mentioning that all EVs are charged to the designated SOCs under the decentralized control. Figures showing the SOC evolutions are omitted here due to the space limit.

\subsection{Network model discrepancies}\label{sec:model_discrepancy}

\begin{figure}[!htb] \centering
\includegraphics[width=0.5\textwidth, trim = 3mm 2mm 12mm 7mm, clip]{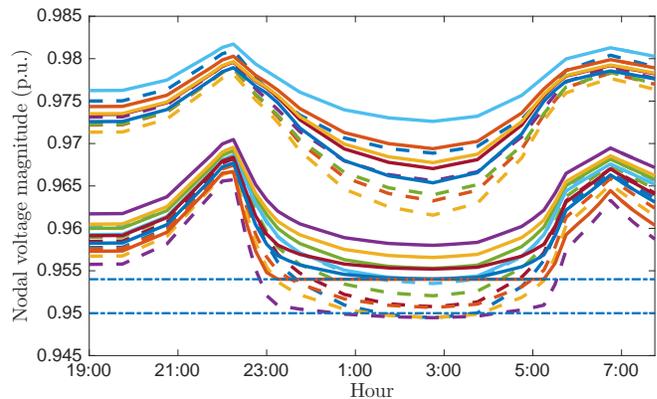}
\caption{Nodal voltage magnitudes calculated from LinDistFlow model (solid lines) and DistFlow model (dashed lines) by adding EV charging loads.}
\label{Nodal_voltage_comp_models}
\end{figure}

Since line losses are discarded in the LinDistFlow model, actual voltage magnitudes that are calculated from the DistFlow model must be lower than those calculated from the LinDistFlow model. Fig. \ref{Nodal_voltage_comp_models} contrasts the voltage differences between LinDistFlow and DistFlow by applying the same decentralized control sequences, where solid and dashed lines represent nodal voltages of LinDistFlow and DistFlow models, respectively.
 It can be observed that, before charging starts at 22:00, voltage discrepancy is as low as 0.001 p.u., while during the charging peak time, the discrepancy is enlarged to about 0.004 p.u.. When nodal voltage magnitudes in LinDistFlow are constrained above 0.954 p.u., lowest magnitudes in DistFlow drop to about 0.95 p.u., which satisfies the ANSI C84.1 standard for service voltage. This suggests that, though the voltage drop serves as constraints in the controller design, considering practical applications, we ought to keep the nodal voltage as close to the ideal value as possible or at least reserve a buffer zone to accommodate line losses.

{\bf{Remark 3:}} In the above simulations, elapsed time of chargers performing one iteration is about 0.25 second. Hence, the control sequences are obtained in 6.5 seconds with 25 iterations. Meanwhile, the centralized algorithm takes 334 seconds. Total computing time and time complexity of the centralized algorithm increases as the number of EVs increases; while due to the scalability, time complexities of \eqref{primal_updates} and \eqref{dual_updates} are independent of EV population size. Specifically, \eqref{primal_updates} is a Euclidean projection onto the intersection of a hyperplane and a hyperbox; \eqref{dual_updates} is a Euclidean projection onto the $\ell_2$-norm ball. These two projections can be solved by a variety of algorithms, among which the most straightforward approach is via quadratic programming. Considering the valley-filling period (7:00 pm - 8:00 am) with 15-min resolution, the complexity of solving a general quadratic programming implementation of \eqref{primal_updates}, vector size $K=52$, is low. Problem size of \eqref{dual_updates} increases with the dimension of the distribution network, thus the complexity might become an issue when a large distribution network is considered. However, solving a projection onto a $\ell_2$-norm ball can be of low complexity if appropriate algorithms are adopted. For example, the Efficient Euclidean Projection in \cite{Su_arxiv_2012} can solve \eqref{dual_updates} with the complexity as low as $O(hK)$, where $h=12$ and $K=52$ in our simulations. Furthermore, runtime issues caused by extremely long control horizon can be resolved by following a model predictive control fashion; runtime issues caused by extremely large feeder size can be resolved by splitting the entire feeder into several parts. \hfill $\blacksquare$


{\bf{Remark 4:}} Though the initial intention is to embed the algorithm into charging points, the proposed decentralized algorithm can also be realized in a centralized controller via parallel computing. \hfill $\blacksquare$

{\bf{Remark 5:}} To ensure the robustness of {\bf{Algorithm \ref{Decentralized_Alg}}}, a waiting time can be set by the central operator. Any $\mathcal{U}_i^{(\ell)}$ not received within the waiting time will be treated as packet loss and then be replaced by the value in its last iteration. Preliminary simulation results show that this packet loss does not impact the algorithm convergency. Rigorous proofs will be provided in the future work. \hfill $\blacksquare$

\section{Concluding Remarks} \label{Conclusion}
This paper developed a decentralized EV charging control framework for the provision of valley-filling in the context of residential distribution network. We formulated the control problem as an optimization problem consisting of a non-separable objective function subject to local constraints and strongly coupled linear inequality network constraints. A novel SPDS algorithm was proposed to solve the formulated problem in a decentralized way. The SPDS ensures convergency of primal and dual variables without regularizing the Lagrangian; the decentralized control scheme allows all chargers update their computations in a parallel fashion and no communication network is needed among chargers. We verified the proposed SPDS algorithm and the decentralized charging control scheme with simulation results.

The proposed SPDS-based decentralized control scheme is not necessarily limited to solely controlling EVs for valley-filling, instead, it is also compatible with DERs and reactive power supplies, and is readily applicable for other grid-level services, e.g., minimization of energy cost. Extending the developed control framework by considering more energy sources and for other grid-level services will be the first step of our future work. Secondly, other constraints such as the transformer overloading could be added to this approach. Unlike the nodal voltages, transformer overloading constraints would be locally coupled ones. This will be discussed in our future work. In addition, we would like to extend the results in this work to a stochastic MPC fashion which can handle the stochastic process of EVs' arrivals and departures. We also intend to study a nonlinear decentralized control scheme based on the DistFlow equations.  Finally, the information required to execute either the centralized or decentralized control scheme includes an accurate network model and knowledge of injections and extractions of real and reactive power at every point in the network.  In practice this information will always be imprecise, and real time control strategies to make up for these errors will be necessary; in the future we will explore the interaction between real-time controllers (e.g., \cite{arnold2016model, baker2017network}) and the scheduling algorithm developed in this paper.

\appendix[Proof of Theorem 2]
Let $\zeta^*=\left[{\mathcal{U}^*}^{\mathsf{T}}~{\lambda^*}^{\mathsf{T}}\right]^\mathsf{T}$ denote the optimizer. By using 
the decomposable structure of $\mathbb{U}$ and the non-expansive property of $\Pi_{\mathbb{X}}(\bold{x})$, we have
\begin{equation} \label{U_relation}
\begin{aligned}
&\left\| \mathcal{U}^{(\ell+1)}-\mathcal{U}^*\right\|_2^2 \\
\leq&\left\| \frac{1}{\tau_\mathcal{U}} \Pi_{\mathbb{U}}\left( \tau_\mathcal{U}\mathcal{U}^{(\ell)}-\alpha\nabla_{\mathcal{U}}\mathcal{L}(\mathcal{U}^{(\ell)},\lambda^{(\ell)})\right) \right. \\
& ~~~~~~ - \frac{1}{\tau_\mathcal{U}} \Pi_{\mathbb{U}}\left( \tau_\mathcal{U}\mathcal{U}^*-\alpha\nabla_{\mathcal{U}}\mathcal{L}(\mathcal{U}^*,\lambda^*)\right) \Big\|_2^2 \\
\leq & \frac{1}{\tau_\mathcal{U}^2}\left\|\tau_\mathcal{U}\mathcal{U}^{(\ell)}-\alpha\nabla_{\mathcal{U}}\mathcal{L}(\mathcal{U}^{(\ell)},\lambda^{(\ell)})-\tau_\mathcal{U}\mathcal{U}_i^*+\alpha\nabla_{\mathcal{U}}\mathcal{L}(\mathcal{U}^*,\lambda^*) \right\|_2^2 \\
=&\frac{1}{\tau_\mathcal{U}^2}\left\| \mathcal{U}^{(\ell)}-\mathcal{U}^* \right\|_2^2 + \frac{\alpha^2}{\tau_\mathcal{U}^2} \left\| \Phi_1(\zeta^{(\ell)})-\Phi_1(\zeta^*) \right\|_2^2 \\
&-2\frac{\alpha}{\tau_\mathcal{U}^2}\left( \Phi_1(\zeta^{(\ell)})-\Phi_1(\zeta^*) \right)^{\mathsf{T}}(\mathcal{U}^{(\ell)}-\mathcal{U}^{*}).
\end{aligned}
\end{equation}

Similarly, we have
\begin{equation} \label{lambda_relation}
\begin{aligned}
&\left\| \lambda^{(\ell+1)}-\lambda^*\right\|_2^2 \\
\leq&\left\| \frac{1}{\tau_\lambda} \Pi_{\mathbb{D}}\left( \tau_\lambda\lambda^{(\ell)}+\beta\nabla_{\lambda}\mathcal{L}(\mathcal{U}^{(\ell)},\lambda^{(\ell)})\right) \right. \\
&~~~~~- \frac{1}{\tau_\lambda}\Pi_{\mathbb{D}}\left( \tau_\lambda\lambda^*+\beta\nabla_{\lambda}\mathcal{L}(\mathcal{U}^*,\lambda^*)\right) \Big\|_2^2 \\
\leq &\frac{1}{\tau_\lambda^2} \left\|\tau_\lambda\lambda^{(\ell)}+\beta\nabla_{\lambda}\mathcal{L}(\mathcal{U}^{(\ell)},\lambda^{(\ell)})-\tau_\lambda\lambda^*-\beta\nabla_{\lambda}\mathcal{L}(\mathcal{U}^*,\lambda^*) \right\|_2^2 \\
=& \frac{1}{\tau_\lambda^2}\left\| \lambda^{(\ell)}-\lambda^* \right\|_2^2 + \frac{\beta^2}{\tau_\lambda^2} \left\| \Phi_2(\zeta^{(\ell)})-\Phi_2(\zeta^*) \right\|_2^2 \\
&-2\frac{\beta}{\tau_\lambda^2}\left( \Phi_2(\zeta^{(\ell)})-\Phi_2(\zeta^*) \right)^{\mathsf{T}}(\lambda^{(\ell)}-\lambda^{*}).
\end{aligned}
\end{equation}

Summing \eqref{U_relation} and \eqref{lambda_relation}, we can readily have
\begin{equation} 
\begin{aligned}
&\left\| \zeta^{(\ell+1)}-\zeta^* \right\|_2^2 \\
\leq& \frac{1}{\tau_\mathcal{U}^2}\left\| \mathcal{U}^{(\ell)}-\mathcal{U}^* \right\|_2^2+\frac{1}{\tau_\lambda^2}\left\| \lambda^{(\ell)}-\lambda^* \right\|_2^2 \\
&+\frac{\alpha^2}{\tau_\mathcal{U}^2} \left\| \Phi_1(\zeta^{(\ell)})-\Phi_1(\zeta^*) \right\|_2^2 + \frac{\beta^2}{\tau_\lambda^2} \left\| \Phi_2(\zeta^{(\ell)})-\Phi_2(\zeta^*) \right\|_2^2 \\
&-2\frac{\alpha}{\tau_\mathcal{U}^2}\left( \Phi_1(\zeta^{(\ell)})-\Phi_1(\zeta^*) \right)^{\mathsf{T}}(\mathcal{U}^{(\ell)}-\mathcal{U}^{*})\\
&-2\frac{\beta}{\tau_\lambda^2}\left( \Phi_2(\zeta^{(\ell)})-\Phi_2(\zeta^*) \right)^{\mathsf{T}}(\lambda^{(\ell)}-\lambda^{*}) \nonumber
\end{aligned}
\end{equation}
\begin{equation} \label{base_inequality}
\begin{aligned}
\leq& \max\left\{ \frac{1}{\tau_\mathcal{U}^2},\frac{1}{\tau_\lambda^2} \right\}\left\| \zeta^{(\ell)}-\zeta^* \right\|_2^2 ~~~~~~~~~~~~~~~~~~~~~~~~~~~~~~~~~~\\
&+\max\left\{\frac{\alpha^2}{\tau_\mathcal{U}^2},\frac{\beta^2}{\tau_\lambda^2}\right\}\left\| \Phi(\zeta^{(\ell)})-\Phi(\zeta^*)\right\|_2^2 \\
&-2\frac{\alpha}{\tau_\mathcal{U}^2}\left( \Phi_1(\zeta^{(\ell)})-\Phi_1(\zeta^*) \right)^{\mathsf{T}}(\mathcal{U}^{(\ell)}-\mathcal{U}^{*})\\
&-2\frac{\beta}{\tau_\lambda^2}\left( \Phi_2(\zeta^{(\ell)})-\Phi_2(\zeta^*) \right)^{\mathsf{T}}(\lambda^{(\ell)}-\lambda^{*}).
\end{aligned}
\end{equation}

The proof will be completed by studying two cases, i.e., ${\alpha}/{\tau_\mathcal{U}^2}>{\beta}/{\tau_\lambda^2}$ and ${\alpha}/{\tau_\mathcal{U}^2}<{\beta}/{\tau_\lambda^2}$.

{\bf{Case 1:}} ${\alpha}/{\tau_\mathcal{U}^2}>{\beta}/{\tau_\lambda^2}$

To proceed, subtracting and adding
\begin{equation}
2\frac{\beta}{\tau_\lambda^2}\left( \Phi_1(\zeta^{(\ell)})-\Phi_1(\zeta^*) \right)^{\mathsf{T}}(\mathcal{U}^{(\ell)}-\mathcal{U}^{*})
\end{equation}
to the right hand side of \eqref{base_inequality} yields
\begin{equation} \label{inequality_2nd}
\begin{aligned}
&\left\| \zeta^{(\ell+1)}-\zeta^* \right\|_2^2 \\
 \leq & \max\left\{ \frac{1}{\tau_\mathcal{U}^2},\frac{1}{\tau_\lambda^2} \right\}\left\| \zeta^{(\ell)}-\zeta^* \right\|_2^2 +\frac{\alpha^2}{\tau_\mathcal{U}^2}\left\| \Phi(\zeta^{(\ell)})-\Phi(\zeta^*)\right\|_2^2 \\
 &-2\frac{\beta}{\tau_\lambda^2} \left( \Phi(\zeta^{(\ell)})-\Phi(\zeta^*)\right)^{\mathsf{T}}(\zeta^{(\ell)}-\zeta^*) \\
 & -2\left(\frac{\alpha}{\tau_\mathcal{U}^2}-\frac{\beta}{\tau_\lambda^2}\right)\left( \Phi_1(\zeta^{(\ell)})-\Phi_1(\zeta^*) \right)^{\mathsf{T}}(\mathcal{U}^{(\ell)}-\mathcal{U}^{*}).
\end{aligned}
\end{equation}

We first deal with the third term on the right hand side. Let $d(\mathcal{U})$, defined in \eqref{Centralized_Lag}, be partitioned by
\begin{equation}
d(\mathcal{U})=\left[ d_1(\mathcal{U})~d_2(\mathcal{U})~\cdots~d_{hK}(\mathcal{U})\right]^{\mathsf{T}} \in \mathbb{R}^{hK},
\end{equation}
where
\begin{equation}
d_j(\mathcal{U})=\mathcal{Y}_{b,j}-D_{d,j}^{\mathsf{T}}\mathcal{U},~j=1,\ldots,h, \nonumber
\end{equation}
and $D_{d,j}$ is the $j$th row of the matrix $D_d$. It can be readily obtained that
\begin{equation}
\begin{aligned}
&\left( \Phi(\zeta^{(\ell)})-\Phi(\zeta^*)\right)^{\mathsf{T}}(\zeta^{(\ell)}-\zeta^*) \\
=&\left( \nabla_{\mathcal{U}}\mathcal{G}(\mathcal{U}^{(\ell)})-\nabla_{\mathcal{U}}\mathcal{G}(\mathcal{U}^{*}) \right)^{\mathsf{T}}(\mathcal{U}^{(\ell)}-\mathcal{U}^*) \\
&+\left(\rho+\frac{1-\tau_\mathcal{U}}{\alpha} \right)\left\|\mathcal{U}^{(\ell)}-\mathcal{U}^* \right\|_2^2 \\
&+\sum_{j=1}^{hK}\left( \lambda^{(\ell)}_j \nabla_{\mathcal{U}}d_j(\mathcal{U}^{(\ell)})-\lambda^{*}_j \nabla_{\mathcal{U}}d_j(\mathcal{U}^{*}) \right)^{\mathsf{T}}(\mathcal{U}^{(\ell)}-\mathcal{U}^*) \\
&-(d(\mathcal{U}^{(\ell)})-d(\mathcal{U}^*))^{\mathsf{T}}(\lambda^{(\ell)}-\lambda^*)+\frac{1-\tau_\lambda}{\beta}\left\| \lambda^{(\ell)}-\lambda^{*}\right\|_2^2.
\end{aligned} 
\end{equation}

Since the function $\mathcal{G}(\mathcal{U})$ is convex, it yields
\begin{equation}
\left( \nabla_{\mathcal{U}}\mathcal{G}(\mathcal{U}^{(\ell)})-\nabla_{\mathcal{U}}\mathcal{G}(\mathcal{U}^{*}) \right)^{\mathsf{T}}(\mathcal{U}^{(\ell)}-\mathcal{U}^*)\geq 0.
\end{equation}
Consequently, we have
\begin{equation} \label{deri_a_b_1st}
\begin{aligned}
&\left( \Phi(\zeta^{(\ell)})-\Phi(\zeta^*)\right)^{\mathsf{T}}(\zeta^{(\ell)}-\zeta^*) \\
\geq&\left(\rho+\frac{1-\tau_\mathcal{U}}{\alpha}\right)\left\|\mathcal{U}^{(\ell)}-\mathcal{U}^* \right\|_2^2+\frac{1-\tau_\lambda}{\beta}\left\| \lambda^{(\ell)}-\lambda^{*}\right\|_2^2 \\
&+\sum_{j=1}^{hK}\left( \lambda^{(\ell)}_j \nabla_{\mathcal{U}}d_j(\mathcal{U}^{(\ell)})-\lambda^{*}_j \nabla_{\mathcal{U}}d_j(\mathcal{U}^{*}) \right)^{\mathsf{T}}(\mathcal{U}^{(\ell)}-\mathcal{U}^*) \\
&-(d(\mathcal{U}^{(\ell)})-d(\mathcal{U}^*))^{\mathsf{T}}(\lambda^{(\ell)}-\lambda^*)  \\
=&\left(\rho+\frac{1-\tau_\mathcal{U}}{\alpha}\right)\left\|\mathcal{U}^{(\ell)}-\mathcal{U}^* \right\|_2^2+\frac{1-\tau_\lambda}{\beta}\left\| \lambda^{(\ell)}-\lambda^{*}\right\|_2^2 \\
&+\sum_{j=1}^{hK}\lambda^{(\ell)}_j\left(d_j(\mathcal{U}^*)-d_j(\mathcal{U}^{(\ell)})- {\nabla_{\mathcal{U}}d_j(\mathcal{U}^{(\ell)})}^{\mathsf{T}}(\mathcal{U}^{*}-\mathcal{U}^{(\ell)} \right) \\
&+\sum_{j=1}^{hK}\lambda^{*}_j\left(d_j(\mathcal{U}^{(\ell)})-d_j(\mathcal{U}^{*})- {\nabla_{\mathcal{U}}d_j(\mathcal{U}^*)}^{\mathsf{T}}(\mathcal{U}^{(\ell)}-\mathcal{U}^{*} \right).
\end{aligned}
\end{equation}

Since $d_j(\mathcal{U}),~j=1,\ldots,h$ is convex, it holds
\begin{equation}
d_j(\mathcal{U})-d_j(\mathcal{V})- {\nabla d_j(\mathcal{V})}^{\mathsf{T}}(\mathcal{U}-\mathcal{V}) \geq 0.
\end{equation}
Hence, the last two terms on the right hand side of \eqref{deri_a_b_1st} are both positive, implying that
\begin{equation} \label{Phi_monotone}
\begin{aligned}
&\left( \Phi(\zeta^{(\ell)})-\Phi(\zeta^*)\right)^{\mathsf{T}}(\zeta^{(\ell)}-\zeta^*) \\
\geq & \left(\rho+\frac{1-\tau_\mathcal{U}}{\alpha}\right)\left\|\mathcal{U}^{(\ell)}-\mathcal{U}^* \right\|_2^2+\frac{1-\tau_\lambda}{\beta}\left\| \lambda^{(\ell)}-\lambda^{*}\right\|_2^2 \\
\geq & \min\left\{\rho+\frac{1-\tau_\mathcal{U}}{\alpha}, \frac{1-\tau_\lambda}{\beta} \right\}\left\| \zeta^{(\ell)}-\zeta^* \right\|_2^2 \\
=& c\left\| \zeta^{(\ell)}-\zeta^* \right\|_2^2.
\end{aligned}
\end{equation}
Without the loss of generality, $\zeta^{(\ell)}$ and $\zeta^{*}$ can be replaced by arbitrary $\zeta_1,\zeta_2 \in \mathbb{U}\times\mathbb{D}$, indicating that the mapping $\Phi(\zeta)$ is strongly monotone with the constant
\begin{equation} \label{monotone_constant}
c=\min\{\rho+\frac{1-\tau_\mathcal{U}}{\alpha},\frac{1-\tau_\lambda}{\beta}\}. 
\end{equation}

To evaluate the second term on the right hand side of \eqref{inequality_2nd}, we have
\begin{equation} \label{inequality_2nd_3rd}
\begin{aligned}
&\left\| \Phi(\zeta^{(\ell)})-\Phi(\zeta^*)\right\|_2 \\
= &\left\| \left[ \begin{array}{c}
\Phi_1(\zeta^{\ell})-\Phi_1(\zeta^*) \\
\Phi_2(\zeta^{\ell})-\Phi_2(\zeta^*)
\end{array}\right] \right\|_2 \\
\leq & \left\|  \nabla_{\mathcal{U}}\mathcal{G}(\mathcal{U}^{(\ell)})-\nabla_{\mathcal{U}}\mathcal{G}(\mathcal{U}^{*})  \right\|_2 \\
&+\left\| \nabla_{\mathcal{U}}({\lambda^{(\ell)}}^{\mathsf{T}}d(\mathcal{U}^{(\ell)}))-\nabla_{\mathcal{U}}({\lambda^{*}}^{\mathsf{T}}d(\mathcal{U}^*))  \right\|_2 \\
&+\left( \rho+\frac{1-\tau_\mathcal{U}}{\alpha}\right)\left\| \mathcal{U}^{(\ell)}-\mathcal{U}^*\right\|_2+\left\| d(\mathcal{U}^{(\ell)})-d(\mathcal{U}^*)\right\|_2 \\
&+\frac{1-\tau_\lambda}{\beta} \left\| \lambda^{(\ell)}-\lambda^* \right\|_2.
\end{aligned}
\end{equation}

The first term on the right hand side of \eqref{inequality_2nd_3rd} gives
\begin{equation}
\begin{aligned}
&\left\|  \nabla_{\mathcal{U}}\mathcal{G}(\mathcal{U}^{(\ell)})-\nabla_{\mathcal{U}}\mathcal{G}(\mathcal{U}^{*})  \right\|_2 \\
=& \left\| \tilde{P}^{\mathsf{T}}\tilde{P}(\mathcal{U}^{(\ell)}-\mathcal{U}^*)\right\|_2\\
\leq & L_{\nabla\mathcal{G}}\left\| \mathcal{U}^{(\ell)}-\mathcal{U}^* \right\|_2,
\end{aligned}
\end{equation}
where $\tilde{P}$ is defined in \eqref{Centralized_objective} and $L_{\nabla\mathcal{G}}=nK\max_{i=1,\ldots,h}\{\bar{P}_i^2 \}$. Arbitrarily choosing $\mathcal{U}_1$ and $\mathcal{U}_2$ to replace $\mathcal{U}^{(\ell)}$ and $\mathcal{U}^*$, we end up with the Lipschitz continuity of $ \nabla_{\mathcal{U}}\mathcal{G}(\mathcal{U})$ with the Lipschitz constant $L_{\nabla\mathcal{G}}$. 

The next step follows the mean-value theorem of vector-valued functions as follows.

{\bf{Theorem 3:}} \emph{Mean-Value Theorem \cite{Apostol_book}: Let $S \subseteq \mathbb{R}^n$ and the mapping $f:~S\mapsto \mathbb{R}^m$ is differentiable at each point of $S$. Let $\bold{x}$ and $\bold{y}$ be two points in $S$ such that all points between $\bold{x}$ and $\bold{y}$ are in $S$. Then for every vector $\bold{a}\in \mathbb{R}^m$, there is a point $\bold{z}$ between $\bold{x}$ and $\bold{y}$ such that
\begin{equation}
\bold{a}^{\mathsf{T}}(f(\bold{x})-f(\bold{y}))=\bold{a}^{\mathsf{T}}(\nabla f(\bold{z})(f(\bold{x})-f(\bold{y}))).
\end{equation}
Further, if $\bold{a}$ is a unit vector such that $\left\| \bold{a}\right\|_2=1$, then it holds that
\begin{equation}
\begin{aligned}
\left\| f(\bold{x})-f(\bold{y} ) \right\|_2 & \leq \left\| \nabla f(\bold{z})(f(\bold{x})-f(\bold{y})) \right\|_2 \\
& \leq L_f \left\| f(\bold{x})-f(\bold{y}) \right\|_2,
\end{aligned}
\end{equation} 
where 
\begin{equation}
\sum_{j=1}^m\left\| \nabla f_j(\bold{z}) \right\|_2 \leq L_f.
\end{equation}}
\hfill $\Box$

Applying {\bf{Theorem 3}} to the mapping $d(\mathcal{U}):\mathbb{R}^{nK}\mapsto \mathbb{R}^{hK}$ yields
\begin{equation} \label{d_Lipschitz}
\left\| d(\mathcal{U}^{(\ell)})-d(\mathcal{U}^*)\right\|_2 \leq L_d \left\| \mathcal{U}^{(\ell)}-\mathcal{U}^*\right\|_2, 
\end{equation}
where $\sum_{j=1}^{hK}\left\| \nabla_{\mathcal{U}} d_j(\mathcal{U}) \right\|_2 \leq L_d$. Since $\nabla_{\mathcal{U}} d_j(\mathcal{U})=D_{d,j}^{\mathsf{T}}$, we have
\begin{equation}
\left\| \nabla_{\mathcal{U}} d_j(\mathcal{U}) \right\|_2 \leq \max_{j=1,\ldots,hK}\left\|D_{d,j}^{\mathsf{T}} \right\|_2,
\end{equation}
and
\begin{equation}
\sum_{j=1}^{hK}\left\| \nabla_{\mathcal{U}} d_j(\mathcal{U}) \right\|_2 \leq hK\max_{j=1,\ldots,hK}\left\|D_{d,j}^{\mathsf{T}} \right\|_2 =L_d.
\end{equation}
Further we can readily have
\begin{equation}
\sum_{j=1}^{hK}\left\| \nabla_{\mathcal{U}} d_j(\mathcal{U}) \right\|_2^2 \leq  \left(\sum_{j=1}^{hK}\left\| \nabla_{\mathcal{U}} d_j(\mathcal{U}) \right\|_2\right)^2 \leq L_d^2.
\end{equation}
This indicates that the mapping $d(\mathcal{U}):\mathbb{R}^{nK}\mapsto \mathbb{R}^{hK}$ is Lipschitz with the constant $L_d$.

Consequently, \eqref{inequality_2nd_3rd} becomes
\begin{equation} \label{Phi_Lipschitz}
\begin{aligned}
&\left\| \Phi(\zeta^{(\ell)})-\Phi(\zeta^*)\right\|_2 \\
\leq & \left(\rho+\frac{1-\tau_\mathcal{U}}{\alpha}+L_{\nabla\mathcal{G}}+L_d\right)\left\| \mathcal{U}^{(\ell)}-\mathcal{U}^*\right\|_2 \\
&+\frac{1-\tau_\lambda}{\beta} \left\| \lambda^{(\ell)}-\lambda^* \right\|_2 \\
&+\left\| \nabla_{\mathcal{U}}({\lambda^{(\ell)}}^{\mathsf{T}}d(\mathcal{U}^{(\ell)}))-\nabla_{\mathcal{U}}({\lambda^{*}}^{\mathsf{T}}d(\mathcal{U}^*))  \right\|_2 \\
=&L_\mathcal{U}\left\| \mathcal{U}^{(\ell)}-\mathcal{U}^*\right\|_2 +\frac{1-\tau_\lambda}{\beta} \left\| \lambda^{(\ell)}-\lambda^* \right\|_2 \\
&+\left\| \sum_{j=1}^{hK}(\lambda^{(\ell)}-\lambda^*)D_{d,j}^{\mathsf{T}} \right\|_2 \\
\leq &L_\mathcal{U}\left\| \mathcal{U}^{(\ell)}-\mathcal{U}^*\right\|_2 +\frac{1-\tau_\lambda}{\beta} \left\| \lambda^{(\ell)}-\lambda^* \right\|_2 \\
&+\sum_{j=1}^{hK}\left|\lambda^{(\ell)}-\lambda^* \right| \left\|D_{d,j}^{\mathsf{T}} \right\|_2 \\
\leq & L_\mathcal{U}\left\| \mathcal{U}^{(\ell)}-\mathcal{U}^*\right\|_2 +\frac{1-\tau_\lambda}{\beta} \left\| \lambda^{(\ell)}-\lambda^* \right\|_2 \\
&+L_d\left\|\lambda^{(\ell)}-\lambda^* \right\|_2 \\
=& L_\mathcal{U}\left\| \mathcal{U}^{(\ell)}-\mathcal{U}^*\right\|_2 +L_\lambda\left\| \lambda^{(\ell)}-\lambda^* \right\|_2 \\
\leq & L_{\Phi}\left\| \zeta^{(\ell)}-\zeta^*\right\|_2,
\end{aligned}
\end{equation}
where
\begin{equation}
L_{\Phi}=\left\|[L_\mathcal{U}, L_\lambda] \right\|_2,
\end{equation}
and
\begin{equation}
\begin{aligned}
L_\mathcal{U}&=\rho+\frac{1-\tau_\mathcal{U}}{\alpha}+L_{\nabla\mathcal{G}}+L_d, \\
L_\lambda&=\frac{1-\tau_\lambda}{\beta}+L_d. 
\end{aligned}
\end{equation}

Next, dropping the coefficient, the last term on the right hand side of \eqref{inequality_2nd} becomes
\begin{equation} \label{inequality_3rd_2nd}
\begin{aligned}
&\left( \Phi_1(\zeta^{(\ell)})-\Phi_1(\zeta^*) \right)^{\mathsf{T}}(\mathcal{U}^{(\ell)}-\mathcal{U}^{*}) \\
=&\left( \nabla_{\mathcal{U}}\mathcal{L}(\mathcal{U}^{(\ell)},\lambda^{(\ell)})-\nabla_{\mathcal{U}}\mathcal{L}(\mathcal{U}^*,\lambda^{(\ell)}) \right)^{\mathsf{T}}\left( \mathcal{U}^{(\ell)}-\mathcal{U}^*\right) \\
&+\left( \nabla_{\mathcal{U}}\mathcal{L}(\mathcal{U}^{*},\lambda^{(\ell)})-\nabla_{\mathcal{U}}\mathcal{L}(\mathcal{U}^*,\lambda^*) \right)^{\mathsf{T}}\left( \mathcal{U}^{(\ell)}-\mathcal{U}^*\right) \\
&+\frac{1-\tau_{\mathcal{U}}}{\alpha}\left\| \mathcal{U}^{(\ell)}-\mathcal{U}^* \right\|_2^2.
\end{aligned}
\end{equation}

The first term on the right hand side of \eqref{inequality_3rd_2nd} indicates that
\begin{equation} \label{inequality_convex_same_D}
\begin{aligned}
&\left( \nabla_{\mathcal{U}}\mathcal{L}(\mathcal{U}^{(\ell)},\lambda^{(\ell)})-\nabla_{\mathcal{U}}\mathcal{L}(\mathcal{U}^*,\lambda^{(\ell)}) \right)^{\mathsf{T}}\left( \mathcal{U}^{(\ell)}-\mathcal{U}^*\right) \\
=&\left( \nabla_{\mathcal{U}}\mathcal{G}(\mathcal{U}^{(\ell)})-\nabla_{\mathcal{U}}\mathcal{G}(\mathcal{U}^*) \right)^{\mathsf{T}}(\mathcal{U}^{(\ell)}-\mathcal{U}^*) +\rho \left\| \mathcal{U}^{(\ell)}-\mathcal{U}^*\right\|_2^2 \\
\geq & \rho \left\| \mathcal{U}^{(\ell)}-\mathcal{U}^*\right\|_2^2.
\end{aligned}
\end{equation}
The equation in \eqref{inequality_convex_same_D} comes from the fact that
\begin{equation}
\nabla_{\mathcal{U}}\left( \lambda^{(\ell)}d( \mathcal{U}^{(\ell)})\right)=\nabla_{\mathcal{U}}\left( \lambda^{(\ell)}d( \mathcal{U}^*)\right).
\end{equation}

The second term on the right hand side of \eqref{inequality_3rd_2nd} indicates that
\begin{equation}  \label{inequality_4th}
\begin{aligned}
&\left( \nabla_{\mathcal{U}}\mathcal{L}(\mathcal{U}^{*},\lambda^{(\ell)})-\nabla_{\mathcal{U}}\mathcal{L}(\mathcal{U}^*,\lambda^*) \right)^{\mathsf{T}}\left( \mathcal{U}^{(\ell)}-\mathcal{U}^*\right) \\
=&\left(  \nabla_{\mathcal{U}}\left({\lambda^{(\ell)}}^{\mathsf{T}}d(\mathcal{U}^*) \right) -   \nabla_{\mathcal{U}}\left({\lambda^{*}}^{\mathsf{T}}d(\mathcal{U}^*) \right) \right)^{\mathsf{T}}\left( \mathcal{U}^{(\ell)}-\mathcal{U}^*\right) \\
=& \left( \sum_{j=1}^{hK}\left({\lambda_j^{(\ell)}}-{\lambda_j^{*}}\right)^{\mathsf{T}}\nabla_{\mathcal{U}}d_j(\mathcal{U}^*) \right)^{\mathsf{T}}  \left( \mathcal{U}^{(\ell)}-\mathcal{U}^*\right) \\
\geq & -\frac{1}{2}\left( \left\| \sum_{j=1}^{hK}\left({\lambda_j^{(\ell)}}-{\lambda_j^{*}}\right)D_{d,j}^{\mathsf{T}} \right\|_2^2 +\left\| \mathcal{U}^{(\ell)}-\mathcal{U}^*\right\|_2^2 \right) \\
\geq &  -\frac{1}{2}\left(  \left( \sum_{j=1}^{hK} \left\| \left({\lambda_j^{(\ell)}}-{\lambda_j^{*}}\right)D_{d,j}^{\mathsf{T}} \right\|_2\right)^2 +\left\| \mathcal{U}^{(\ell)}-\mathcal{U}^*\right\|_2^2 \right) \\
= &  -\frac{1}{2}\left(  \left( \sum_{j=1}^{hK} \left\| {\lambda_j^{(\ell)}}-{\lambda_j^{*}}\right\|_2 \left\|D_{d,j}^{\mathsf{T}} \right\|_2\right)^2  +\left\| \mathcal{U}^{(\ell)}-\mathcal{U}^*\right\|_2^2 \right) \\
= & -\frac{1}{2}  \left(\left[ \begin{array}{c}
\left\| \lambda_1^{(\ell)}-\lambda_1^*\right\|_2\\
\left\| \lambda_2^{(\ell)}-\lambda_2^*\right\|_2\\
\vdots \\
\left\| \lambda_{hK}^{(\ell)}-\lambda_{hK}^*\right\|_2
\end{array} \right]^\mathsf{T} \left[  \begin{array}{c}
\left\| D_{d,1}^{\mathsf{T}} \right\|_2 \\
\left\| D_{d,2}^{\mathsf{T}} \right\|_2 \\
\vdots \\
\left\| D_{d,h}^{\mathsf{T}} \right\|_2
\end{array} \right] \right)^2  \\
&~ -\frac{1}{2}\left\| \mathcal{U}^{(\ell)}-\mathcal{U}^*\right\|_2^2  \\
\geq & -\frac{1}{2}  \left( \left\| \lambda^{(\ell)}-\lambda^* \right\|_2 \left\|\left[  \begin{array}{c}
\left\| D_{d,1}^{\mathsf{T}} \right\|_2 \\
\left\| D_{d,2}^{\mathsf{T}} \right\|_2 \\
\vdots \\
\left\| D_{d,h}^{\mathsf{T}} \right\|_2
\end{array} \right] \right\|_2 \right)^2 \\
&~  -\frac{1}{2}\left\| \mathcal{U}^{(\ell)}-\mathcal{U}^*\right\|_2^2 \\
=& -\frac{1}{2}\left( \sum_{j=1}^{hK}\left\| D_{d,j}^{\mathsf{T}} \right\|_2^2 \left\| \lambda^{(\ell)}-\lambda^* \right\|_2^2 + \left\| \mathcal{U}^{(\ell)}-\mathcal{U}^*\right\|_2^2 \right) \\
\geq & -\frac{1}{2}\left( L_d^2\left\| \lambda^{(\ell)}-\lambda^* \right\|_2^2 + \left\| \mathcal{U}^{(\ell)}-\mathcal{U}^*\right\|_2^2  \right). 
\end{aligned}
\end{equation}

Substituting \eqref{Phi_monotone}, \eqref{Phi_Lipschitz}, \eqref{inequality_convex_same_D}, and \eqref{inequality_4th} into \eqref{inequality_2nd}, we have
\begin{equation} 
\begin{aligned}
&\left\| \zeta^{(\ell+1)}-\zeta^* \right\|_2^2 \\
 \leq &\left(\max\left\{\frac{1}{\tau_{\mathcal{U}}^2},\frac{1}{\tau_\lambda^2} \right\}+\frac{\alpha^2}{\tau_{\mathcal{U}}^2}L_\Phi^2-2c\frac{\beta}{\tau_\lambda^2} \right. \\
 &\left.+\left(\frac{\alpha}{\tau_\mathcal{U}^2}-\frac{\beta}{\tau_\lambda^2} \right)\max\left\{L_d^2,1-2\left(\rho+\frac{1-\tau_\mathcal{U}}{\alpha}\right)\right\} \right) \\
 &\times \left\| \zeta^{(\ell)}-\zeta^* \right\|_2^2 \\
 =& \varrho  \left\| \zeta^{(\ell)}-\zeta^* \right\|_2^2.
\end{aligned}
\end{equation}
This coincides with \eqref{zeta_convergency}.

{\bf{Case 2:}} ${\alpha}/{\tau_\mathcal{U}^2}<{\beta}/{\tau_\lambda^2}$

Inequality \eqref{base_inequality} can be rewritten as
\begin{equation} \label{inequality_5th}
\begin{aligned}
&\left\| \zeta^{(\ell+1)}-\zeta^* \right\|_2^2 \\
\leq &\max\left\{\frac{1}{\tau_\mathcal{U}^2},\frac{1}{\tau_\lambda^2} \right\} \left\| \zeta^{(\ell)}-\zeta^* \right\|_2^2 +\frac{\beta^2}{\tau_\lambda^2}\left\| \Phi(\zeta^{(\ell)})-\Phi(\zeta^*)\right\|_2^2 \\
 &-2\frac{\alpha}{\tau_\mathcal{U}^2} \left( \Phi(\zeta^{(\ell)})-\Phi(\zeta^*)\right)^{\mathsf{T}}(\zeta^{(\ell)}-\zeta^*) \\
 &+2(\frac{\alpha}{\tau_\mathcal{U}^2}-\frac{\beta}{\tau_\lambda^2})\left( \Phi_2(\zeta^{(\ell)})-\Phi_2(\zeta^*) \right)^{\mathsf{T}}(\lambda^{(\ell)}-\lambda^{*}) \\
 \leq &\left(\max\left\{\frac{1}{\tau_\mathcal{U}^2},\frac{1}{\tau_\lambda^2} \right\}+\frac{\beta^2}{\tau_\lambda^2}L_{\Phi}^2-2c\frac{\alpha}{\tau_\mathcal{U}^2} \right) \left\| \zeta^{(\ell)}-\zeta^* \right\|_2^2 \\
 & +2(\frac{\beta}{\tau_\lambda^2}-\frac{\alpha}{\tau_\mathcal{U}^2})\left( \nabla_{\lambda}\mathcal{L}(\mathcal{U}^{(\ell)},\lambda^{(\ell)})-\nabla_{\lambda}\mathcal{L}(\mathcal{U}^*,\lambda^*) \right)^{\mathsf{T}} \\
 &\times \left( \lambda^{(\ell)}-\lambda^*\right) + 2(\frac{\alpha}{\tau_\mathcal{U}^2}-\frac{\beta}{\tau_\lambda^2})\frac{1-\tau_\lambda}{\beta}\left\| \lambda^{(\ell)}-\lambda^* \right\|_2^2.
\end{aligned}
\end{equation}
If we can further prove that the second term on the right hand side of \eqref{inequality_5th} is smaller or equal to $l\left\| \zeta^{(\ell)}-\zeta^*\right\|_2^2$, where $l\geq0$, then the proof is complete. By applying Cauchy-Schwarz inequality, we have
\begin{equation} \label{inequality_6th}
\begin{aligned}
&\left( \nabla_{\lambda}\mathcal{L}(\mathcal{U}^{(\ell)},\lambda^{(\ell)})-\nabla_{\lambda}\mathcal{L}(\mathcal{U}^*,\lambda^*) \right)^{\mathsf{T}}\left( \lambda^{(\ell)}-\lambda^*\right) \\
\leq & \left| \left( \nabla_{\lambda}\mathcal{L}(\mathcal{U}^{(\ell)},\lambda^{(\ell)})-\nabla_{\lambda}\mathcal{L}(\mathcal{U}^*,\lambda^*) \right)^{\mathsf{T}}\left( \lambda^{(\ell)}-\lambda^*\right) \right| \\
\leq & \left\|  \nabla_{\lambda}\mathcal{L}(\mathcal{U}^{(\ell)},\lambda^{(\ell)})-\nabla_{\lambda}\mathcal{L}(\mathcal{U}^*,\lambda^*)  \right\|_2 \left\| \lambda^{(\ell)}-\lambda^* \right\|_2 \\
\leq & \frac{1}{2} \left\|  \nabla_{\lambda}\mathcal{L}(\mathcal{U}^{(\ell)},\lambda^{(\ell)})-\nabla_{\lambda}\mathcal{L}(\mathcal{U}^*,\lambda^*)  \right\|_2^2  + \frac{1}{2}\left\| \lambda^{(\ell)}-\lambda^* \right\|_2^2 \\
=& \frac{1}{2}\left\| d({\mathcal{U}}^{(\ell)})-d(\mathcal{U}^*) \right\|_2^2+ \frac{1}{2}\left\| \lambda^{(\ell)}-\lambda^* \right\|_2^2.
\end{aligned}
\end{equation}

By applying {\bf{Theorem 3}} and the boundness on $\nabla d(\mathcal{U})$ to \eqref{inequality_6th}, we have
\begin{equation} \label{inequality_7th}
\begin{aligned}
&\left( \nabla_{\lambda}\mathcal{L}(\mathcal{U}^{(\ell)},\lambda^{(\ell)})-\nabla_{\lambda}\mathcal{L}(\mathcal{U}^*,\lambda^*) \right)^{\mathsf{T}}\left( \lambda^{(\ell)}-\lambda^*\right) \\
\leq & \frac{1}{2}L_d^2 \left\| {\mathcal{U}}^{(\ell)}-\mathcal{U}^*\right\|_2^2+ \frac{1}{2}\left\| \lambda^{(\ell)}-\lambda^* \right\|_2^2.
\end{aligned}
\end{equation}

Substituting \eqref{inequality_7th} into \eqref{inequality_5th}, we can readily have
\begin{equation}
\begin{aligned}
&\left\| \zeta^{(\ell+1)}-\zeta^* \right\|_2^2 \\
 \leq &\left(\max\left\{\frac{1}{\tau_{\mathcal{U}}^2},\frac{1}{\tau_\lambda^2} \right\}+\frac{\beta^2}{\tau_{\lambda}^2}L_\Phi^2-2c\frac{\alpha}{\tau_\mathcal{U}^2} \right. \\
 &\left.+\left(\frac{\beta}{\tau_\lambda^2}-\frac{\alpha}{\tau_\mathcal{U}^2} \right)\max\left\{L_d^2,1-2\frac{1-\tau_\lambda}{\beta}\right\} \right) \\
 &\times \left\| \zeta^{(\ell)}-\zeta^* \right\|_2^2 \\
 =& \varrho  \left\| \zeta^{(\ell)}-\zeta^* \right\|_2^2.
\end{aligned}
\end{equation}
This coincides with \eqref{zeta_convergency}. 

Further, we prove $\exists~\alpha>0,\beta>0$ such that $0<\varrho<1$. Let $\psi$ denote $\max\{\hat{\alpha}, \hat{\beta} \}$, let $\min\{\hat{\alpha},\hat{\beta}\}=\mu \max\{\hat{\alpha}, \hat{\beta}\}$, where $0 < \mu < 1$, and let $\omega$ denote $\max\left\{ {1}/{\tau_\mathcal{U}^2}, {1}/{\tau_\lambda^2}\right\}$. Then the expression of $\varrho$ becomes
\begin{equation}
\varrho=\omega+\delta\psi L_{\Phi}^2-2c\mu\psi+(1-\mu)\psi \phi.
\end{equation}

In order to make the sequence $\{\zeta^{(\ell)}\}$ converge, we need to have $0<\varrho<1$, implying that
\begin{subequations}
\begin{equation} \label{1st_situation}
\omega+\delta\psi L_{\Phi}^2-2c\mu\psi+(1-\mu)\psi \phi< 1,
\end{equation}
\begin{equation} \label{2nd_situation}
\omega+\delta\psi L_{\Phi}^2-2c\mu\psi+(1-\mu)\psi \phi>0.
\end{equation}
\end{subequations}
Solving \eqref{1st_situation} for $\mu$ gives
\begin{equation}
\mu \in \left(\frac{(\omega-1)+\delta\psi L_\Phi^2+\phi\psi}{2c\psi+\phi\psi},\infty\right) \cap (0,1).
\end{equation}
To have a feasible solution in $(0,1)$, we need to have
\begin{equation} \label{converge_condition}
\omega+\delta\psi L_\Phi^2-1 < 2c\psi.
\end{equation}

Solving \eqref{2nd_situation} for $\mu$ gives
\begin{equation}
\mu \in \left(-\infty, \frac{\omega+\delta\psi L_\Phi^2+\phi\psi}{2c\psi+\phi\psi}\right) \cap (0,1),
\end{equation}
whose solution always exists given that \eqref{converge_condition} holds. Thus, we can always find a tuple $(\alpha, \beta)$ having
\begin{equation}
\mu \in \left(\frac{(\omega-1)+\delta\psi L_\Phi^2+\phi\psi}{2c\psi+\phi\psi}, \min\left\{ \frac{\omega+\delta\psi L_\Phi^2+\phi\psi}{2c\psi+\phi\psi}   ,1  \right\}\right),
\end{equation}
indicating that $0<\varrho<1$. Then the sequence of $\{\zeta^{(\ell)}\}$ converges. This completes the proof. \hfill $\blacksquare$

%





\ifCLASSOPTIONcaptionsoff
  \newpage
\fi



\bibliographystyle{IEEEtran}
\end{document}

%% file: sy-headinput.tex
\newcommand{\ba}{\begin{array}}
\newcommand{\ea}{\end{array}}
\newcommand{\bc}{\begin{center}}
\newcommand{\ec}{\end{center}}

\newcommand{\beqn}[1]{\begin{equation}\label{#1}}
\newcommand{\eeqn}{\end{equation}}
\newcommand{\be}{\begin{equation}}
\newcommand{\ee}{\end{equation}}

\newcommand{\beqnn}{\begin{eqnarray}}
\newcommand{\eeqnn}{\end{eqnarray}}

\newcommand{\sgn}{{\rm sgn}}

